\documentclass[reqno,11pt]{amsart}
\usepackage{amsfonts}

\usepackage{dsfont,mathscinet}

\usepackage[left=2.8cm,right=2.8cm,top=2.8cm,bottom=2.8cm]{geometry}
\usepackage{graphicx}

\usepackage{epstopdf}
\usepackage[colorlinks=true, pdfstartview=FitV, linkcolor=blue, 
            citecolor=blue, urlcolor=blue]{hyperref}

\usepackage[T1]{fontenc}
\usepackage{lmodern}
\usepackage{microtype}
\usepackage{amsthm}
\theoremstyle{plain}
\newtheorem{theorem}{Theorem}[section]
\newtheorem{proposition}[theorem]{Proposition}
\newtheorem{lemma}[theorem]{Lemma}
\newtheorem{coro}[theorem]{Corollary}

\theoremstyle{definition}
\newtheorem{definition}[theorem]{Definition}
\newtheorem{remark}[theorem]{Remark}
\newtheorem*{remark*}{Remark}

\newtheorem*{hypo*}{Hypothesis}

\usepackage{amsmath,amssymb,mathrsfs,mathtools,braket,latexsym,esint}
\numberwithin{equation}{section}

\def\nat{\mathbb{N}}
\def\inte{\mathbb{Z}}

\def\real{\mathbb{R}}

\def\matr{\real^{3 \times 3}}

\def\cD{\mathcal{D}}
\def\cE{\mathcal{E}}
\def\cF{\mathcal{F}}
\def\cG{\mathcal{G}}
\def\cJ{\mathcal{J}}

\def\rmQ{\mathcal{Q}'}
\def\rmT{\mathrm{T}}

\def\ny{\nabla y}
\def\nv{\nabla v}

\def\nP{\nabla P}
\def\el{\mathrm{el}}

\def\homo{\mathrm{hom}}
\def\per{\mathrm{per}}
\def\SL{\mathsf{SL}}

\def\tr{\mathrm{tr}}
\def\io{\int_\Omega}

\def\dd{\,\mathrm{d}}
\def\wts{\stackrel{2}{\rightharpoonup}}

\def\ep{\varepsilon}
\def\R{\mathbb{R}}

\def\eps{\varepsilon}
\renewcommand{\epsilon}{\eps}

\usepackage{enumerate}
\usepackage[shortlabels]{enumitem}

\usepackage{color}

\usepackage{tikz}

\title[Homogenization of high-contrast media]{Homogenization of high-contrast media\\in finite-strain elastoplasticity}
\date{}

\author[E. Davoli]{Elisa Davoli}
\author[C. Gavioli]{Chiara Gavioli}
\author[V. Pagliari]{Valerio Pagliari}	

\AtEndDocument{
	\bigskip
	\footnotesize{
		\noindent
		E. Davoli, C. Gavioli, V. Pagliari.
		
		\noindent
		Institute of Analysis and Scientific Computing,
		TU Wien,
		
		\noindent
		Wiedner Hauptstra\ss e 8-10, 1040 Vienna, Austria.
		
		\noindent
		E-mails:
		\href{mailto:elisa.davoli@tuwien.ac.at}{\tt elisa.davoli@tuwien.ac.at},
		\href{mailto:chiara.gavioli@tuwien.ac.at}{\tt chiara.gavioli@tuwien.ac.at}, \href{mailto:valerio.pagliari@tuwien.ac.at}{\tt valerio.pagliari@tuwien.ac.at}.
	}
}

\begin{document}
	
	\begin{abstract}
	    This work is devoted to the analysis of the interplay
	    between internal variables and high-contrast microstructure
	    in inelastic solids. As a concrete case-study,
	    by means of variational techniques,
		we derive a macroscopic description
		for an elastoplastic medium.
		Specifically, we consider a composite
		obtained by filling the voids of a periodically perforated stiff matrix
		by soft inclusions.
		We study the $\Gamma$-convergence of the related energy functionals
		as the periodicity tends to zero,
		the main challenge being posed by the lack of coercivity
		brought about by the degeneracy of the material properties
		in the soft part.
		We prove that the $\Gamma$-limit,
		which we compute with respect to a suitable notion of convergence,
		is the sum of the contributions
		resulting from each of the two components separately.
		Eventually, convergence of the energy minimizing configurations is obtained.
		
		\medskip
		\noindent
		{\it 2020 Mathematics Subject Classification:}
		49J45; 
		74B20; 
		74C15; 
		74E30; 
		74Q05. 
		
		\smallskip
		\noindent
		{\it Keywords and phrases:}
		finite-strain elastoplasticity, $\Gamma$-convergence, homogenization, high-contrast, two-scale convergence.
	\end{abstract}
	
	\maketitle
	
	{\parskip=0em \tableofcontents}

	\section{Introduction}\label{sec:intro}
	The present paper is concerned
	with the variational analysis of some integral functionals
	that model the stored energy
	of materials governed by finite-strain elastoplasticity with hardening.
	Our goal is to derive, by means of $\Gamma$-convergence,
	the effective macroscopic energy
	of a special class of heterogeneous materials,
	those with a so called high-contrast microstructure.
	The interest in such media stems from the experimental observation of an infinite number of band gaps in their mechanical behavior:
	high-contrast materials, indeed, exhibit infinitely many interval of frequencies in which wave propagation is not allowed.
	This, in turn, makes them extremely interesting for possible cloaking applications.
	Some recent ones in civil engineering, for example in seismic waves cloaking,
	and in the modeling of advanced sensor and actuator devices
	call for advancements in the mathematical modeling
	of those classes of high-contrast materials
	that have not been fully studied yet,
	like the ones we consider here.
	
	The mathematical literature on high-contrast materials is vast. To keep our presentation concise, we only point out that, besides results for stratified elastoplastic composites \cite{cc,cc1,cdf,cd}, the only additional available contributions in the inelastic setting concern the study of brittle fracture problems \cite{barchiesi, barchiesi.lazzaroni.zeppieri, pellet.scardia.zeppieri}. For the modeling of nonlinear elastic high-contrast composites we single out the works \cite{braides.garroni, CC}.
	
	When undertaking the analysis of high-contrast media beyond the elastic purview, 
	hurdles are posed by the mathematical treatment of possible internal variables and dissipative effects, as well as by their interplay with the high-contrast microstructure. In this paper we initiate such task by focusing on the case-study of finite elastoplasticity (see, e.g., \cite{lubliner}). At this first stage we neglect both the difficulties due to possible lack of coercivity for the dissipative effects and those associated with time evolution.
	Thus, we focus here on a static model for a single time-step with a global regularization on the gradient of the plastic strain, and leave the analysis of different regimes and the passage to the limit in the quasistatic evolutions for future investigations.
	
	The present study grounds
	on a previous result that we obtained in \cite{DGP1},
	where we addressed the static homogenization
	of elastoplastic microstructures in the large strain regime.	
	As in that work, our starting point is the description of the medium at the microscopic level.
	We let $\Omega \subset \real^3$ be an open, bounded, connected set with Lipschitz boundary,
	and we suppose it to be the reference configuration
	of an elastoplastic body that exhibits the following microstructure:
	denoting by $\eps>0$ the microscale,
	we suppose that a stiff perforated matrix $\Omega^1_\eps$ sits in $\Omega$
	and that its pores are filled by soft inclusions,
	which form the set $\Omega^0_\eps$ (see Figure~\ref{fig:Omega}).
	Let us denote by $\SL(3)$
	the group of $3 \times 3$ real matrices with determinant equal to~$1$.
	When the matrix and the inclusions exhibit the same plastic-hardening $H$,	
	the functionals encoding the stored energy associated with
	the deformation $y\in W^{1,2}(\Omega;\real^3)$ and the plastic strain $P \in W^{1,q}(\Omega;K)$, with $q>3$ and $K\subset\SL(3)$ a given compact set, read
	\begin{equation}\label{eq:Jeps}
		\begin{aligned}
		\cJ_\eps (y,P) & \coloneqq
		\int_{\Omega^0_\eps} W^0_\eps\left(\eps\ny(x)P^{-1}(x)\right) \dd x
		+ \int_{\Omega^1_\eps} W^1\left(\ny(x)P^{-1}(x)\right) \dd x \\
		& \quad + \io H\big( P(x) \big ) \dd x+\int_\Omega |\nabla P(x)|^q \dd x,
		\end{aligned}
	\end{equation}
	where $\{W^0_\eps\}_{\eps > 0}$ and $W^1$ are, respectively,
	the elastic energy densities of the inclusions and of the matrix.
	
	Let us briefly comment on some modeling choices
	underlying \eqref{eq:Jeps}.
	The factor $\eps$ multiplying the argument of $W^0_\eps$ encodes
	the high-contrast between the two components, and
	it results in a loss of coercivity in the problem.
	From a modeling perspective,
	this heuristically means that very large deformations of the inclusions are allowed or, 
	in other words, that the inclusions are very soft
	-- whence the expression {\em high-contrast} to describe the difference between the phases.
		
	As for the hardening term, note that
	also additional hardening variables have been taken into account in the literature,
	see \cite{Mie02,mielke} for a modeling overview.
	Here, for the purpose of putting the high-contrast behavior to the foreground,
	we give up full generality and 
	restrict ourselves to the case in which
	only a hardening dependence on the plastic strain is given.
	A discussion on alternative modeling choices is also presented in Remark~\ref{hardening}.
		
	Our main result describes the asymptotics of the functionals $\cJ_\epsilon$,
	and it is presented in Theorem~\ref{stm:Glim}.
	The precise mathematical framework of our analysis is described in Section~\ref{sec:math},
	where further details on the definitions
	and on the roles of the terms in $\cJ_\eps$ may be found.
		
	We work under the classical assumption that the elastic behavior of our sample $\Omega$ is independent of preexistent plastic distortions.
	Then, the deformation gradient $\nabla y$
	associated with a deformation $y \colon \Omega\to \R^3$
	of the body decomposes into an elastic strain and a plastic one.
	In the framework of linearized elastoplasticity
	the decomposition would  take an additive form.
	In the case at stake,
	that of finite plasticity \cite{kroner, Lee, mielke, Mie02}, 
	the existence of an intermediate configuration
	determined by purely plastic distortions
	is instead assumed, and
	it is then supposed that
	elastic deformations are applied to such intermediate configuration.
	Mathematically, these hypotheses amount to 
	a multiplicative decomposition of the gradient of a deformation $y\in W^{1,2}(\Omega;\R^3)$:
	$$\ny(x) = F_\el(x)P(x) \quad \text{for a.\,e. }x \in \Omega,$$
	for a suitable \textit{elastic strain} $F_\el \in L^2(\Omega;\matr)$ and a \textit{plastic strain}
	$P \in L^2(\Omega;{\rm SL}(3))$.
	On the one hand, such multiplicative structure has recently found an atomistic validation
	in the framework of crystal plasticity
	by means of a discrete-to-continuum analysis
	\cite{conti.reina, conti.schloemerkamper.reina}.	
	On the other hand, 
	alternative models for finite plasticity have been proposed.
	However, since a discussion of fine modeling issues goes beyond the scopes of our work,
	we do not dwell here on a comparison of the various modeling theories.
	We refer the reader interested in this topic
	to, e.g., \cite{davoli.francfort, grandi.stefanelli1, grandi.stefanelli2, naghdi}.

	Finally, we comment on the regularizing term in $\nabla P$ in the energy \eqref{eq:Jeps}. As mentioned before, at this stage we assume it to provide coercivity of the energy with respect to the plastic-strain variables on the whole set $\Omega$.
	From a modeling point of view, we note that this regularization is common in engineering models, for it prevents the formation of microstructures, see \cite{BCHH,CHM}. Alternative higher order regularizations are discussed in \cite{drs}.
	

	Let us conclude our introduction with a few words on the proofs.
	A delicate point is choosing a suitable notion of convergence that ensures effective compactness properties. Indeed, the fact that the energy contribution in the soft inclusions is evaluated in terms of $\ep \nabla y$ leads to a loss of coercivity, and, subsequently, to the loss of compactness in classical weak Sobolev topologies. On the other hand, using the strong two-scale convergence of the gradients (as in \cite{CC}) does not guarantee convergence of minimizers of $\cJ_\eps$ to minimizers of the limiting functional. To cope with this difficulty, we adapt the approach of \cite{DKP} and introduce an ad hoc notion of convergence for deformations, to which we refer as \emph{convergence in the sense of extensions}. Roughly speaking, a sequence of deformations converges in the sense of extensions if it is bounded in $L^2$ and can be extended in $W^{1,2}$ in such a way that the extensions are weakly compact in the Sobolev sense, cf.~Definition~\ref{stm:defconv} and Remarks~\ref{stm:conv-tilde} and \ref{stm:traces} for the precise definition and some basic properties. For the plastic strains, we argue instead by using the uniform convergence.
	This choice is motivated by the fact that sequences of deformations and plastic strains with uniformly bounded energies are precompact
	with respect to the convergence resulting from pairing the two mentioned above.
	Thus our $\Gamma$-convergence analysis directly entails convergence of minimizers.
	We observe that this result can easily be extended to functionals which take into account also plastic dissipation. We refer to Section~\ref{sec:conclusions} for a more detailed discussion on this point.
	
	Our approach to the proofs resorts to extension results on perforated domains, to two-scale convergence and periodic unfolding techniques, as well as to equiintegrability arguments used to control the behavior of the microstructure close to the boundary of the set $\Omega$.
	A key step is a splitting procedure that allows to treat the soft and the stiff parts separately.
	
	\subsection*{Outline of the paper}
	The setup of our analysis and the main result, Theorem~\ref{stm:Glim},
	are presented in Section~\ref{sec:math}.
	Section~\ref{sec:pre} contains some useful preliminaries.
	In Section~\ref{sec:comp+split} we discuss the equicoercivity
	of the energy functionals under consideration and the splitting procedure.
	The asymptotic behavior of the soft inclusions is characterized in Section~\ref{sec:Glimsoft}.
	The ground is then laid for the proof of Theorem~\ref{stm:Glim}, which is contained in Section~\ref{sec:conclusions} together with a variant including plastic dissipation and a comparison
	with the aforementioned result from \cite{CC}.
	
	\section{Mathematical setting and results}\label{sec:math}
	
	Hereafter, $\Omega$ is an open, bounded, and connected set
	with Lipschitz boundary in $\real^3$.
	Working in the $3$-dimensional space is not essential,
	and our analysis can be easily adapted
	to the setting of $\real^d$ with $d=2$ or $d>3$.
	The spaces of real-valued $3\times 3$ and $3\times 3 \times 3$ tensors
	are denoted by $\matr$ and  by $\real^{3\times 3 \times 3}$, respectively.
	We adopt the symbol $I$ for the identity matrix.
	By $|\,\cdot\,|$ we denote indiscriminately
	the Euclidean norms in  $\real^3$, $\matr$ and $\real^{3\times 3 \times 3}$.
	To deal with plastic strains, we recall the classical notation
		\begin{align*}
			\SL(3) \coloneqq \{F \in \real^{3 \times 3} : \det F = 1\}. 
		\end{align*}
	
	If $A \subset \real^3$ is a measurable set,
	we denote by $\mathcal{L}^3(A)$ its three-dimensional Lebesgue measure.
	Finally, we simply write $\|\,\cdot\,\|_{L^p}$ for the norm of a function
	in $L^p(\Omega;\real^3)$, $L^p(\Omega;\real^{3\times 3})$, or $L^p(\Omega;\real^{3\times 3\times 3})$
	when no ambiguity arises,
	and we specify the integration domain only when necessary.
	
	A fundamental role in  our study is played
	by the following notion of variational convergence,
	see the monograph \cite{DalM} for a thorough treatment:
	\begin{definition}\label{def:Gammaconv}
		Let $X$ be a set endowed with a notion of convergence. We say that a family of functionals $\{\cG_\eps\}$,
		with $\cG_\eps \colon X \to [-\infty,+\infty]$, \textit{$\Gamma$-converges} as $\eps \to 0$ to $\cG\colon X \to [-\infty,+\infty]$ if for all $x \in X$ and all infinitesimal sequences $\{\eps_k\}_{k\in\nat}$ the following holds:
		\begin{enumerate}
			\item for every sequence $\{x_k\}_{k\in\nat} \subset X$ such that $x_k \to x$, we have
			$$\cG(x) \le \liminf_{k\to+\infty} \cG_{\eps_k}(x_k);$$
			\item there exists a sequence $\{x_k\}_{k\in\nat} \subset X$ such that $x_k \to x$ and
			$$\limsup_{k\to+\infty} \cG_{\eps_k}(x_k) \le \cG(x).$$
		\end{enumerate}
	\end{definition}
	When $X$ is equipped with a topology $\tau$,
	we write e.g.\ $\Gamma(\tau)$-convergence
	to stress what the underlying convergence for sequences in $X$ is.
	In what follows, for notational convenience,
	we indicate the dependence on $\eps_k$ by means of the subscript $k$ alone,
	e.g., $\cJ_k \coloneqq \cJ_{\eps_k}$.
	
	Our aim is to study elastoplastic media with high-contrast periodic microstructure
	in the case of soft inclusions inserted in a perforated stiff matrix.
	Letting $ Q\coloneqq (0,1)^3$ be the periodicity cell,
	in order to describe the geometry in precise terms,
	we start by considering	a set $E^1\subset \R^3$ 
	that is open, connected, $Q$-periodic, and has Lipschitz boundary,
	cf.~\cite{ACDP} or~\cite[Chapter 19]{BrDFr}.
	We recall that the set $E^1 \subset \real^3$  is $Q$-periodic if
	$E^1+t=E^1$ for all $t\in \inte^3$.
	The set $E^1$ is then employed to define the microstructure as follows.
	
\begin{figure}
		\caption{The periodicity cell $Q$ and
			its partition into the soft inclusion $Q^0$ (white) 
			and the stiff matrix $Q^1$ (gray).
		}
		\label{fig:periodcell}		
		\centering
		\bigskip
		
		\begin{tikzpicture}[scale=0.75]
			\filldraw[fill=lightgray] (-2,-2) rectangle (2,2);
			
			\filldraw[fill=white,xscale=2]
			(-0.7,0.9) .. controls (-0.6,1) and (-0.4,1.4) ..
			(-0.2,1.3) .. controls (0,1.2) and (0.1,1.1) ..
			(0.2,0.6)  .. controls (0.3,0.2) and (-0.2,-0.6) .. 
			(-0.5,0.1) .. controls (-0.8,0.6) and (-0.8,0.8) ..
			(-0.7,0.9);
			
			\filldraw[fill=white,rotate around={-120:(1.1,-0.1)},scale=1.5]
			(1.4,-0.3) .. controls (1,0.2) and (1.1,-0.3) ..
			(0.9,0) .. controls (0.7,0.3) and (0.9,0.4) ..
			(1.4,0.3) .. controls (1.9,0.1) and (1.8,-0.8) ..
			(1.4,-0.3);
			
			\path (-0.5,0.6) node {$Q^0$};
			\path (0.8,-1) node {$Q^0$};
			\path (1.2,1) node {$Q^1$};
		\end{tikzpicture}
	\end{figure}
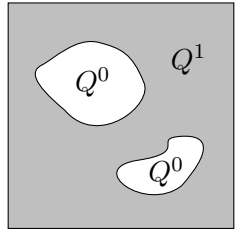

	First, at the scale level $\eps = 1$, we define
	\[
	Q^1\coloneqq Q\cap E^1
	\quad\text{and}\quad
	Q^0\coloneqq Q\setminus  \overline{Q^1},
	\]
	where the sets $Q^0$ and $Q^1$ represent, respectively,
	the inclusion and the matrix of the unit cell $Q$
	(see Figure~\ref{fig:periodcell}).
	Note that, according to the definition of $Q^1$, it holds that
	\begin{equation}\label{eq:E1}
		E^1 \coloneqq \bigcup_{t \in \mathbb{Z}^3} (t+Q^1).
	\end{equation}
	
\begin{figure}
		\caption{The microstructure of the composite in $\Omega$.
			The soft inclusions that form $\Omega^0_\eps$ correspond to the white holes,
			while the grey region represents the matrix $\Omega^1_\eps$.
            Note that the perforations do not intersect the boundary.}
		\label{fig:Omega}
		\centering
		
		\bigskip	
		\begin{tikzpicture}
			\filldraw[thick,fill=lightgray,rotate=30,scale=1.5]
			[xshift=4,yshift=-2] (-1.5,1)  .. controls (-1.5,1.5) and (-1,2) ..
			(0,2.5) .. controls (1,3) and (1,3) ..
			(1.5,3) .. controls (2.25,3) and (3.5,1.5) .. 
			(3,1) .. controls (2.5,0.5) and (2,1) .. 
			(2.5,0)  .. controls (3,-1) and (0.5,-1) ..
			(0,-0.5) .. controls (-0.5,0) and (-1.5,0.5) ..
			(-1.5,1);
			
			\draw[very thin,<->,] (-3.2,-1) -- (-3.2,0);
			\path (-3.4,-0.5) node {$\varepsilon$};
			
			\draw[very thin] (-3,-1) -- (3,-1);
			\foreach \y in {0,1,...,3}
			{
				\draw[very thin] (-3,\y) -- (4,\y);
			}
			\draw[very thin] (-2,4) -- (4,4);
			\draw[very thin] (-2,5) -- (4,5);
			\draw[very thin] (-1,6) -- (2,6);
			
			\draw[very thin] (-3,-1) -- (-3,3);
			\draw[very thin] (-2,-1) -- (-2,5);
			\foreach \x in {-1,0,...,2}
			{
				\draw[very thin] (\x,-1) -- (\x,6);
			}
			\draw[very thin] (3,-1) -- (3,5);
			\draw[very thin] (4,0) -- (4,5);
			
			\foreach \x in {-1.5,-0.5,...,2.5}
			\foreach \y in {0.5,1.5}
			{
				\filldraw[very thin,fill=white,xshift=\x cm,yshift=\y cm,xscale=0.5,yscale=0.25]
				(-0.7,0.9) .. controls (-0.6,1) and (-0.4,1.4) ..
				(-0.2,1.3) .. controls (0,1.2) and (0.1,1.1) ..
				(0.2,0.6)  .. controls (0.3,0.2) and (-0.2,-0.6) .. 
				(-0.5,0.1) .. controls (-0.8,0.6) and (-0.8,0.8) ..
				(-0.7,0.9);
				
				\filldraw[very thin,fill=white,xshift=\x cm,yshift=\y cm,scale=0.25,rotate around={-120:(1.1,0.1)},scale=1.5]
				(1.4,-0.3) .. controls (1,0.2) and (1.1,-0.3) ..
				(0.9,0) .. controls (0.7,0.3) and (0.9,0.4) ..
				(1.4,0.3) .. controls (1.9,0.1) and (1.8,-0.8) ..
				(1.4,-0.3);
			}
			
			\foreach \x in {-1.5,-0.5,...,1.5}
			\foreach \y in {2.5}
			{
				\filldraw[very thin,fill=white,xshift=\x cm,yshift=\y cm,xscale=0.5,yscale=0.25]
				(-0.7,0.9) .. controls (-0.6,1) and (-0.4,1.4) ..
				(-0.2,1.3) .. controls (0,1.2) and (0.1,1.1) ..
				(0.2,0.6)  .. controls (0.3,0.2) and (-0.2,-0.6) .. 
				(-0.5,0.1) .. controls (-0.8,0.6) and (-0.8,0.8) ..
				(-0.7,0.9);
				
				\filldraw[very thin,fill=white,xshift=\x cm,yshift=\y cm,scale=0.25,rotate around={-120:(1.1,0.1)},scale=1.5]
				(1.4,-0.3) .. controls (1,0.2) and (1.1,-0.3) ..
				(0.9,0) .. controls (0.7,0.3) and (0.9,0.4) ..
				(1.4,0.3) .. controls (1.9,0.1) and (1.8,-0.8) ..
				(1.4,-0.3);
			}
			
			\foreach \x in {-0.5,0.5,...,2.5}
			\foreach \y in {3.5}
			{
				\filldraw[very thin,fill=white,xshift=\x cm,yshift=\y cm,xscale=0.5,yscale=0.25]
				(-0.7,0.9) .. controls (-0.6,1) and (-0.4,1.4) ..
				(-0.2,1.3) .. controls (0,1.2) and (0.1,1.1) ..
				(0.2,0.6)  .. controls (0.3,0.2) and (-0.2,-0.6) .. 
				(-0.5,0.1) .. controls (-0.8,0.6) and (-0.8,0.8) ..
				(-0.7,0.9);
				
				\filldraw[very thin,fill=white,xshift=\x cm,yshift=\y cm,scale=0.25,rotate around={-120:(1.1,0.1)},scale=1.5]
				(1.4,-0.3) .. controls (1,0.2) and (1.1,-0.3) ..
				(0.9,0) .. controls (0.7,0.3) and (0.9,0.4) ..
				(1.4,0.3) .. controls (1.9,0.1) and (1.8,-0.8) ..
				(1.4,-0.3);
			}
			
			\foreach \x in {0.5,1.5}
			\foreach \y in {4.5}
			{
				\filldraw[very thin,fill=white,xshift=\x cm,yshift=\y cm,xscale=0.5,yscale=0.25]
				(-0.7,0.9) .. controls (-0.6,1) and (-0.4,1.4) ..
				(-0.2,1.3) .. controls (0,1.2) and (0.1,1.1) ..
				(0.2,0.6)  .. controls (0.3,0.2) and (-0.2,-0.6) .. 
				(-0.5,0.1) .. controls (-0.8,0.6) and (-0.8,0.8) ..
				(-0.7,0.9);
				
				\filldraw[very thin,fill=white,xshift=\x cm,yshift=\y cm,scale=0.25,rotate around={-120:(1.1,0.1)},scale=1.5]
				(1.4,-0.3) .. controls (1,0.2) and (1.1,-0.3) ..
				(0.9,0) .. controls (0.7,0.3) and (0.9,0.4) ..
				(1.4,0.3) .. controls (1.9,0.1) and (1.8,-0.8) ..
				(1.4,-0.3);
			}
			
			
			\foreach \x in {-2.5,-1.5,...,2.5}
			\foreach \y in {-0.5}
			{
				\draw[very thin,xshift=\x cm,yshift=\y cm,xscale=0.5,yscale=0.25]
				(-0.7,0.9) .. controls (-0.6,1) and (-0.4,1.4) ..
				(-0.2,1.3) .. controls (0,1.2) and (0.1,1.1) ..
				(0.2,0.6)  .. controls (0.3,0.2) and (-0.2,-0.6) .. 
				(-0.5,0.1) .. controls (-0.8,0.6) and (-0.8,0.8) ..
				(-0.7,0.9);
				
				\draw[very thin,xshift=\x cm,yshift=\y cm,scale=0.25,rotate around={-120:(1.1,0.1)},scale=1.5]
				(1.4,-0.3) .. controls (1,0.2) and (1.1,-0.3) ..
				(0.9,0) .. controls (0.7,0.3) and (0.9,0.4) ..
				(1.4,0.3) .. controls (1.9,0.1) and (1.8,-0.8) ..
				(1.4,-0.3);
			}
			
			\foreach \x in {-2.5,3.5}
			\foreach \y in {0.5,1.5,2.5}
			{
				\draw[very thin,xshift=\x cm,yshift=\y cm,xscale=0.5,yscale=0.25]
				(-0.7,0.9) .. controls (-0.6,1) and (-0.4,1.4) ..
				(-0.2,1.3) .. controls (0,1.2) and (0.1,1.1) ..
				(0.2,0.6)  .. controls (0.3,0.2) and (-0.2,-0.6) .. 
				(-0.5,0.1) .. controls (-0.8,0.6) and (-0.8,0.8) ..
				(-0.7,0.9);
				
				\draw[very thin,xshift=\x cm,yshift=\y cm,scale=0.25,rotate around={-120:(1.1,0.1)},scale=1.5]
				(1.4,-0.3) .. controls (1,0.2) and (1.1,-0.3) ..
				(0.9,0) .. controls (0.7,0.3) and (0.9,0.4) ..
				(1.4,0.3) .. controls (1.9,0.1) and (1.8,-0.8) ..
				(1.4,-0.3);
			}
			
			{
				\draw[very thin,xshift=2.5 cm,yshift=2.5 cm,xscale=0.5,yscale=0.25]
				(-0.7,0.9) .. controls (-0.6,1) and (-0.4,1.4) ..
				(-0.2,1.3) .. controls (0,1.2) and (0.1,1.1) ..
				(0.2,0.6)  .. controls (0.3,0.2) and (-0.2,-0.6) .. 
				(-0.5,0.1) .. controls (-0.8,0.6) and (-0.8,0.8) ..
				(-0.7,0.9);
				
				\draw[very thin,xshift=2.5 cm,yshift=2.5 cm,scale=0.25,rotate around={-120:(1.1,0.1)},scale=1.5]
				(1.4,-0.3) .. controls (1,0.2) and (1.1,-0.3) ..
				(0.9,0) .. controls (0.7,0.3) and (0.9,0.4) ..
				(1.4,0.3) .. controls (1.9,0.1) and (1.8,-0.8) ..
				(1.4,-0.3);
			}
			
			\foreach \x in {-1.5,3.5}
			\foreach \y in {3.5,4.5}
			{
				\draw[very thin,xshift=\x cm,yshift=\y cm,xscale=0.5,yscale=0.25]
				(-0.7,0.9) .. controls (-0.6,1) and (-0.4,1.4) ..
				(-0.2,1.3) .. controls (0,1.2) and (0.1,1.1) ..
				(0.2,0.6)  .. controls (0.3,0.2) and (-0.2,-0.6) .. 
				(-0.5,0.1) .. controls (-0.8,0.6) and (-0.8,0.8) ..
				(-0.7,0.9);
				
				\draw[very thin,xshift=\x cm,yshift=\y cm,scale=0.25,rotate around={-120:(1.1,0.1)},scale=1.5]
				(1.4,-0.3) .. controls (1,0.2) and (1.1,-0.3) ..
				(0.9,0) .. controls (0.7,0.3) and (0.9,0.4) ..
				(1.4,0.3) .. controls (1.9,0.1) and (1.8,-0.8) ..
				(1.4,-0.3);
			}	
			
			\foreach \x in {-0.5,2.5}
			\foreach \y in {4.5}
			{
				\draw[very thin,xshift=\x cm,yshift=\y cm,xscale=0.5,yscale=0.25]
				(-0.7,0.9) .. controls (-0.6,1) and (-0.4,1.4) ..
				(-0.2,1.3) .. controls (0,1.2) and (0.1,1.1) ..
				(0.2,0.6)  .. controls (0.3,0.2) and (-0.2,-0.6) .. 
				(-0.5,0.1) .. controls (-0.8,0.6) and (-0.8,0.8) ..
				(-0.7,0.9);
				
				\draw[very thin,xshift=\x cm,yshift=\y cm,scale=0.25,rotate around={-120:(1.1,0.1)},scale=1.5]
				(1.4,-0.3) .. controls (1,0.2) and (1.1,-0.3) ..
				(0.9,0) .. controls (0.7,0.3) and (0.9,0.4) ..
				(1.4,0.3) .. controls (1.9,0.1) and (1.8,-0.8) ..
				(1.4,-0.3);
			}
			
			\foreach \x in {-0.5,0.5,1.5}
			\foreach \y in {5.5}
			{
				\draw[very thin,xshift=\x cm,yshift=\y cm,xscale=0.5,yscale=0.25]
				(-0.7,0.9) .. controls (-0.6,1) and (-0.4,1.4) ..
				(-0.2,1.3) .. controls (0,1.2) and (0.1,1.1) ..
				(0.2,0.6)  .. controls (0.3,0.2) and (-0.2,-0.6) .. 
				(-0.5,0.1) .. controls (-0.8,0.6) and (-0.8,0.8) ..
				(-0.7,0.9);
				
				\draw[very thin,xshift=\x cm,yshift=\y cm,scale=0.25,rotate around={-120:(1.1,0.1)},scale=1.5]
				(1.4,-0.3) .. controls (1,0.2) and (1.1,-0.3) ..
				(0.9,0) .. controls (0.7,0.3) and (0.9,0.4) ..
				(1.4,0.3) .. controls (1.9,0.1) and (1.8,-0.8) ..
				(1.4,-0.3);
			}
		\end{tikzpicture}
	\end{figure}
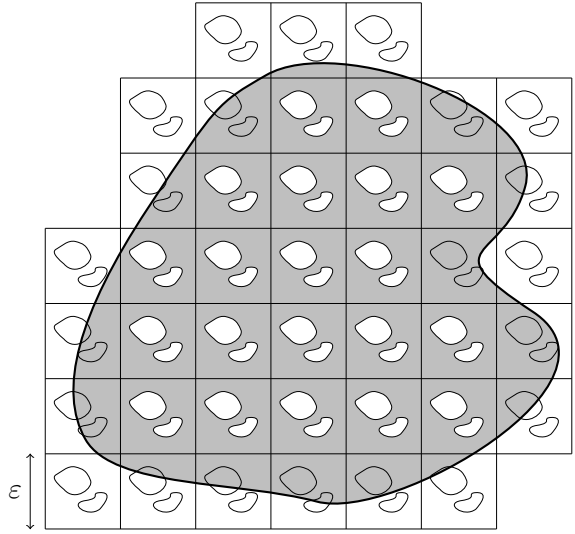

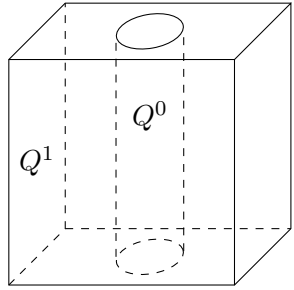
\begin{figure}
		\caption{In the $3$-dimensional space, interconnected soft fibers do not disconnect the matrix.
			A simple case is depicted here:
			the cylindrical perforation $Q^0$ runs through the periodicity cell
			and its complement $Q^1$ is connected.}
		\label{fig:fibers}
		
		\bigskip
		
		\centering
		
		\begin{tikzpicture}[scale=0.75]
			\draw (-2,-2) rectangle (2,2);
			\draw (2,-2) -- (3,-1) -- (3,3) -- (2,2);
			\draw (-2,2) -- (-1,3) -- (3,3);
			\draw[dashed] (-2,-2) -- (-1,-1) -- (-1,3);
			\draw[dashed] (-1,-1) -- (3,-1);
			
			\draw (.5,2.5) ellipse [x radius=.6,y radius=.3,rotate=10]; 
			\draw[dashed] (.5,-1.5) ellipse [x radius=.6,y radius=.3,rotate=10];
			\draw[dashed] (1.1,-1.4) -- (1.1,2.5);
			\draw[dashed] (-0.1,-1.5) -- (-0.1,2.4);
			
			\path (0.5,1) node {$Q^0$};
			\path (-1.5,0.2) node {$Q^1$};
		\end{tikzpicture}
		
	\end{figure}

	Second, $Q^0$ is translated and rescaled to describe the set of soft inclusions.
	Precisely, given a (small) $\lambda>0$,
	we define the collection of inclusions at a scale $\eps>0$ as
	\begin{equation}\label{eq:Omega0}
		\Omega_\eps^0 \coloneqq \bigcup_{t \in T_\eps}  \eps(t+Q^0),
		\quad\text{where }
		T_\eps\coloneqq \Big\{t\in \inte^3 : \mathrm{dist}\Big(\eps (t+Q^0), \partial\Omega \Big) > \lambda \eps \Big\}.
	\end{equation}
	Since $\Omega$ represents the region of space occupied by the whole composite,
	the stiff matrix is then given by
	\begin{gather}
		\label{eq:Omega1}
		\Omega_\eps^1 \coloneqq \Omega \setminus \overline{\Omega_\eps^0},
	\end{gather}
	see Figure~\ref{fig:Omega}.
	Note that the set $\Omega_\eps^1$ is connected and Lipschitz, and that
	\eqref{eq:Omega0} ensures that
	the inclusions are compactly contained in $\Omega$,
	since they are separated from the boundary by a strip of width $\lambda \eps$.
	Our assumptions allow for some flexibility on the geometry of the inclusions,
	which could for instance form interconnected fibers (see Figure~\ref{fig:fibers}).
	Indeed,	differently from other works (e.g.~\cite{CC}),
	we do not prescribe that the unit perforation $Q_0$ is compactly contained in $Q$.
	Therefore, the geometry considered in this paper is on the one hand less restrictive
	than that in the seminal contributions on perforated domains,
	but on the other hand it is less general than that in~\cite{ACDP}.
	
	Our $\Gamma$-convergence result deals with the asymptotic behavior,
	as $\eps$ tends to $0$, of the family $\{\cJ_\eps\}$ defined by \eqref{eq:Jeps}.
	Before stating the result,
	we collect the hypotheses we use in the following lines.
	
	The elastic energy density of the stiff matrix
	$W^1 \colon \matr \to [0,+\infty]$ satisfies the following:
		\begin{enumerate}[label=\textbf{E\arabic*:},ref={E\arabic*}]
			\item\label{E-growth} It is $2$-coercive and has at most quadratic growth, that is,
			there exist $0 < c_1 \le c_2$ such that for all $F \in \matr$
			$$ c_1 |F|^2 \le W^1(F) \le c_2\left(|F|^2+1\right).$$
			\item\label{E-lip} It is $2$-Lipschitz:
			there exists $c_3 > 0$ such that for all $F_1,F_2 \in \matr$
			$$	|W^1(F_1) - W^1(F_2)| \le c_3 \left(1 + |F_1| + |F_2|\right)|F_1 - F_2|. $$
		\end{enumerate}	
	
	The assumptions on the soft energy densities $W^0_\eps \colon \matr \to [0,+\infty]$ are analogous:
	\begin{enumerate}[label=\textbf{E\arabic*:},ref={E\arabic*},start=3]
		\item\label{E-growth-eps} There exist $0 < c_1 \le c_2$ such that for all $F \in \matr$, and all $\eps > 0$,
		\begin{align*}
			c_1 |F|^2 \le W^0_\eps(F) \le c_2\left(|F|^2+1\right).
		\end{align*}
		\item\label{E-lip-eps} There exists $c_3 > 0$ such that for all $F_1,F_2 \in \matr$, and all $\eps > 0$,
		\begin{align*}
			\left|W^0_\eps(F_1)-W^0_\eps(F_2)\right| \le c_3 \left(1+|F_1|+|F_2|\right)|F_1-F_2|.
		\end{align*}
		\item\label{E-conv}  There exists $W^0 \colon \matr \to [0,+\infty]$ such that for all $F \in \matr$
		$$\lim_{\eps \to 0} W^0_\eps(F) = W^0(F).$$
	\end{enumerate}
	
	\begin{remark}\label{regW0}
		The function $W^0$ possesses the same growth and regularity properties of $W^0_\eps$.
	\end{remark}
	Our assumptions rule out non-impenetrability constraints
	at the level of the energy.
	A blow-up of the energy on matrices with non-positive determinant is desirable from a modeling point of view, but at the same time very difficult to handle and yet to be done in the context of homogenization.
	Frame indifference is instead compatible with our hypotheses,
	and up to adding a constant, conditions \ref{E-growth} and \ref{E-growth-eps} fulfill the physical requirement that rigid motions have zero elastic energy.	
	We also note that the choice of
	considering a family $\{W^0_\eps\}$ instead of a fixed $W^0$ for the soft stored elastic energy
	is somehow standard in the literature,
	see, e.g., \cite[Remark~2]{CC} for a motivational example in the context of solid mechanics.
	
	Next, we list the assumptions on the hardening $H \colon \matr \to [0,+\infty]$.
		\begin{enumerate}[label=\textbf{H\arabic*:},ref={H\arabic*}]
			\item\label{H1} Assume that a \textit{Finsler structure} on $\SL(3)$ is assigned.
				$H(F)$ is finite if and only if $F \in K$,
				where $K\subset \SL(3)$ is a geodesically convex, compact neighborhood of $I$.
			\item\label{H2} The restriction of $H$ to $K$ is Lipschitz continuous.
		\end{enumerate}
	The requirement that $K$ is geodesically convex
	with respect to the Finsler structure assigned on $\SL(3)$
	is the crucial ingredient to invoke \cite[Theorem~2.2]{DGP1},
	which in our context is employed to
	capture the asymptotic behavior of the stiff matrix,
	see Theorem~\ref{stm:homo-fin-plast2}.
	We refer to \cite{DGP1} for a discussion on the role
	of the Finsler geometry for the homogenization of elastoplastic media, and
	to Subsection~\ref{sec:Finsler} for a summary of the tools from that theory
	that we need here.
	In particular,
	the existence of a set $K$ complying with \ref{H1} is settled
	in Lemma~\ref{stm:K-in-H2} below.
		
		Requirement \ref{H1} prescribes that
		the effective domain of $H$
		coincides with a compact set $K$ containing $I$.
		Then it follows that there exists $c_K>0$ such that
		\begin{equation}\label{Pbound}
			|F| + |F^{-1}| \le c_K
			\quad \text{for every }F \in K ,
		\end{equation}
		because $\SL(3)$ is by definition well separated from $0$. 
		As a consequence, plastic strains with finite hardening are uniformly bounded in $L^\infty$,
		and, in particular, we infer that
		for any $F \in K$ and $G\in\matr$
		\begin{equation}\label{stm:Pbound}
			\left|G\right| = \left|G F^{-1} F\right| \le c_K\left| G F^{-1}\right|.
		\end{equation}		
	
	\begin{remark}
	\label{hardening}
	Note that in principle it would be reasonable to suppose that
	the soft and the stiff components feature different hardening behaviors.
	For instance, it could be imposed that
	the soft hardening is evaluated on an $\eps$-rescaling of the plastic stress,
	thus replicating the structure of the elastic contribution.
	As the only available tool to deal with periodic homogenization at finite strains is \cite[Theorem~2.2]{DGP1},
	we leave such scenarios for possible future investigation and
	restrict ourselves to a simpler setting,
	namely we choose to model both hardening terms
	by a single function satisfying \ref{H1} and \ref{H2}.
	We point out that under these assumptions
	making a distinction between $H^i = H^i(P)$, $i = 0,1$
	would not require any substantial change in our approach,
	therefore we dispense with it.
	Qualitatively,
	keeping the soft hardening contribution of order $1$ amounts to the situation
	in which, for small $\eps$,
	elastic deformations much larger than the plastic ones are allowed.
	\end{remark}

	We can now state the homogenization result for high-contrast elastoplastic media.
	Since we want our analysis to yield convergence of minima and minimizers of $\cJ_\eps$ to the ones of the limiting energy,
	we need to introduce a convergence that is compliant with the degeneracy of the soft inclusions.
	For shortness, we refer to it as \emph{convergence in the sense of extensions},
	even though the name is not at all standard.
	
	\begin{definition}\label{stm:defconv}
		Let $\{\eps_k\}$ be an infinitesimal sequence.
		We say that $\{y_k\} \subset W^{1,2}(\Omega;\real^3)$ converges to $y\in W^{1,2}(\Omega;\real^3)$
		in the sense of extensions with respect to the scales $\eps_k$
		if the following hold:
		\begin{enumerate}
			\item $\{y_k\}$ is bounded in $ L^2(\Omega;\real^3)$;
			\item there exists a sequence $\{\tilde y_k\}\subset W^{1,2}(\Omega;\real^3)$
			such that $y_k = \tilde y_k$ in $\Omega^1_k\coloneqq \Omega^1_{\eps_k}$ and
			$\tilde y_k \rightharpoonup y$ weakly in $ W^{1,2}(\Omega;\real^3)$.
		\end{enumerate}
	\end{definition}
	
	\begin{remark}\label{stm:conv-tilde}
		Let $\tilde{y}_k = \tilde{y}'_k$ a.e.\ in $\Omega^1_k$.
		Let as well $\tilde{y}_k \to y$
		and $\tilde{y}'_k \to y'$ strongly in $L^2(\Omega;\real^3)$
        (e.g., $y$ and $y'$ are $W^{1,2}(\Omega;\real^3)$-weak limits of the respective sequences).
		Then, recalling \eqref{eq:E1} and \eqref{eq:Omega1}, and observing that $\Omega \cap\eps_k E^1\subset \Omega^1_k$, we get
		$$0 = \lim_{k\to+\infty} \int_{\Omega^1_k} |\tilde{y}_k - \tilde{y}'_k| \dd x \ge \lim_{k\to+\infty} \io \chi_{\eps_k E^1}(x) |\tilde{y}_k - \tilde{y}'_k| \dd x =  \mathcal{L}^3(Q_1)\io |y - y'| \dd x.$$
		From this, we conclude that $y = y'$ a.e.\ in $\Omega$.
		In particular, if the limit in the sense of extensions exists, then it is unique.
	\end{remark}
	
	\begin{remark}\label{stm:traces}
		By \eqref{eq:Omega0}, there exists a neighborhood $O_k$ of $\partial\Omega$ in $\Omega$ such that $\Omega^1_k \cap O_k \equiv \Omega \cap O_k$.
		Therefore, if $y$ and $\tilde{y}$ coincide in $\Omega^1_k$, their traces on $\partial\Omega$ are also equal.
	\end{remark}
	
	Bearing in mind that we set $q>3$,
    the asymptotic behavior of the family $\{\cJ_\eps\}$
    is described in the next theorem:
	
	\begin{theorem}\label{stm:Glim}
		Let $\{W^1\}$ and $\{W^0_\eps\}$ satisfy \ref{E-growth}--\ref{E-conv}, and
		let $H$ satisfy \ref{H1}--\ref{H2}.
		For all $y \in  L^2(\Omega;\real^3)$ and $P \in L^q(\Omega;\SL(3))$ 
		there exists a functional
		$$\cJ(y,P) \coloneqq \Gamma\mbox{-}\lim_{\eps \to 0} \cJ_\eps(y,P),$$
		where the underlying convergences are
		the one in the sense of extensions and the uniform one,
		respectively for the first and for the second argument.
		The $\Gamma$-limit is characterized as follows:
		$$\cJ(y,P) = \cJ^0(P) + \cJ^1(y,P),$$
		\begin{equation}\label{eq:J0}
				\cJ^0(P) \coloneqq \left\lbrace
				\begin{aligned}
					\displaystyle
					\mathcal{L}^3(Q^0) \io \Big[ \rmQ W^0 \big( 0 , P^{-1}(x) \big)  + H\big( P(x) \big ) \Big] \dd x
					&\quad \text{if } 
                        P\in  W^{1,q}(\Omega;K),\\[3pt]
					+\infty \qquad\qquad\qquad\qquad\qquad\qquad\qquad\qquad\quad &\quad \text{otherwise in } 
                    L^q(\Omega;\SL(3)),
				\end{aligned}
				\right.
		\end{equation}
		and
		\begin{equation}\label{eq:J1}
					\cJ^1(y,P) \coloneqq \left\lbrace\begin{aligned}
						\displaystyle\io \Big[\widetilde W^1_\homo\big( \ny(x),P(x) \big) + \mathcal{L}^3(Q^1) H &\big( P(x) \big)  + |\nP(x)|^q\Big] \dd x \\
						& \text{if } (y,P) \in W^{1,2}(\Omega;\real^3)\times W^{1,q}(\Omega;K), \\[3pt]
						+\infty \qquad\qquad\qquad\qquad\qquad\qquad\quad &\text{otherwise in } L^2(\Omega;\real^3)\times L^q(\Omega;\SL(3)).
					\end{aligned}
					\right.			
		\end{equation}
		Here, for $F,G \in \matr$,
		\begin{equation}\label{eq:tildeW0}
			\rmQ W^0(F,G)\coloneqq \inf \left\{ \int_Q 
			W^0 \Big(\big( F + \nv(z) \big) G \Big) \dd z : v\in W^{1,2}_0(Q;\real^3)\right\},
		\end{equation}
		while
		\begin{multline*}
			\widetilde W^1_\homo(F,G) \coloneqq \lim_{\lambda  \to +\infty} \frac{1}{\lambda ^3} \inf\Bigg\{
			\int_{(0,\lambda )^3 \cap E^1} W^1\Big(\big(F+\ny(x)\big)G^{-1}\Big) \dd x 
			: y \in W^{1,2}_0((0,\lambda )^3;\real^3) \Bigg\}.
		\end{multline*}
	\end{theorem}
	
	The formula defining $\rmQ W^0$ provides a variant of the classical quasiconvex envelope of $W^0$.  
	We refer to Section~\ref{sec:Glimsoft} for further discussion on this point.
	
\begin{remark}\label{rk:zero}
In principle, it cannot be excluded 
that some nontrivial energy densities $W^0_\eps$ do not contribute
to the elastic homogenized energy,
in the sense that, when finite,
for the corresponding $\cJ^0$ we have
$$
\cJ^0(P)=\mathcal{L}^3(Q^0) \io H\big( P(x) \big ) \dd x.
$$
As an instance of this phenomenon, 
we consider the following example.
For any $F\in \matr$, we let
$W^0_\eps(F)=W^0(F)\coloneqq |F|^2$.
Conditions~\ref{E-growth-eps}--\ref{E-conv} are satisfied by definition.
Since for any fixed $G\in \matr$
the function $F\mapsto W^0_G(F)\coloneqq W^0(FG)$ is convex, it is, in particular, also quasiconvex. Hence, $\rmQ W^0(0,G)=W^0(0,G)=W^0(0)=0$.
\end{remark}
		
As a byproduct of our asymptotic analysis,
we are in a position to infer convergence
of the minimization problems associated with the energy functionals and
of the related (quasi)~minimizers.
	
	\begin{coro}\label{stm:conv-min}
		Let the assumptions and notation of Theorem~\ref{stm:Glim} hold, and let $\{(y_k,P_k)\} \subset W^{1,2}_0(\Omega;\real^3) \times W^{1,q}(\Omega;K)$ be a sequence of almost minimizers, that is,
		$$\lim_{k\to+\infty} \Big( \cJ_k(y_k,P_k) - \inf \cJ_k(y,P)\Big) = 0,$$
		where the infimum is taken over $W^{1,2}_0(\Omega;\real^3) \times W^{1,q}(\Omega;K)$.
		Then, there exists a minimizer
		$(y,P)\in W^{1,2}_0(\Omega;\real^3) \times W^{1,q}(\Omega;K)$
		of $\cJ$ such that, up to subsequences,
		$y_k \to y$ in the sense of extensions and $P_k \to P$ uniformly. 
		Moreover,
		$$\inf \cJ_k \to \min \cJ.$$
	\end{coro}

	\begin{remark}
		The conclusion of the previous corollary is not affected
		if the homogeneous boundary conditions on $\{y_k\}$ are replaced
		by more general (and physical) ones, for example
		$y_k = u$ on a non negligible subset of $\partial \Omega$ for a given $u \in  W^{1,2}(\Omega;\real^3)$.
	
		Note, instead, that global forcing terms such as
		\[
			\cF^{\rm ext}(y)\coloneqq -\int_\Omega f\cdot y \dd x
			\quad\text{for a given } f\in W^{1,2}(\Omega;\real^3)
		\]
		cannot be added to the functionals $\cJ_k$
		without need of further analysis.
		Indeed, the functional $\cF^{\rm ext}$ is not continuous
		with respect to the convergence in Definition~\ref{stm:defconv}, and
		therefore standard results about continuous perturbations
		in the context of $\Gamma$-convergence cannot be invoked.
	\end{remark}	

The proof of Theorem~\ref{stm:Glim} consists of three steps.
First, we study the compactness properties of sequences $\{(y_\eps,P_\eps)\}$ satisfying $\sup_\eps \cJ_\eps(y_\eps,P_\eps) \le C$ and characterize their limits.
Second, we show that the two components of the material can be studied independently.
Finally, we perform the analysis of each component separately.
In view of this approach, it is useful to introduce the functionals
that account for the two different contributions, namely
	\begin{align}\label{eq:Eeps0}
		\cE_\eps^0 (y,P) & \coloneqq
		\io \chi_\eps^0(x)\Big[
		W^0_\eps\left(\eps\ny(x)P^{-1}(x)\right) + H\big( P(x) \big )
		\Big]\! \dd x, \\
		\label{eq:Eeps1}
		\cE_\eps^1 (y,P) & \coloneqq
		\io \chi_\eps^1(x)\Big[
		W^1\left(\ny(x)P^{-1}(x)\right) + H\big( P(x) \big )
		\Big]\!  \dd x, 
	\end{align}
	where, for $i = 0,1$, $\chi_\eps^i(x)$ denotes the characteristic function of $\Omega_\eps^i$
	(i.e., $\chi_\eps^i(x)=1$ if $x\in  \Omega_\eps^i$ and
	$\chi_\eps^i(x)=0$ otherwise).
	We also decompose the functional $\cJ_\eps$ accordingly:
	\begin{align*}
		\cJ_\eps = \cJ_\eps^0 + \cJ_\eps^1,
	\end{align*}
	with
	\begin{gather}
		\label{eq:Jeps0}
		\cJ_\eps^0(y,P) \coloneqq \begin{cases}
			\cE_\eps^0 (y,P) & \text{if } (y,P) \in W^{1,2}(\Omega;\real^3)\times W^{1,q}(\Omega;K), \\
			+\infty &\text{otherwise in } L^2(\Omega;\real^3)\times L^q(\Omega;\SL(3)),
		\end{cases}
		\\
		\label{eq:Jeps1}
		\cJ_\eps^1(y,P) \coloneqq \begin{cases}
			\cE_\eps^1 (y,P)
			+ \lVert \nabla P \rVert^q_{L^q(\Omega;\real^{3 \times 3 \times 3})} & \text{if } (y,P) \in W^{1,2}(\Omega;\real^3)\times W^{1,q}(\Omega;K), \\
			+\infty &\text{otherwise in } L^2(\Omega;\real^3)\times L^q(\Omega;\SL(3)).
		\end{cases}
	\end{gather}
	
	In contrast to $\cJ_\eps^1(y,P)$,
	whose asymptotic behavior is derived from \cite[Theorem~2.2]{DGP1},
	the soft part requires a dedicated treatment.
	This happens already in the setting of nonlinear elasticity (see \cite{CC}).
We obtain the following proposition, whose proof is given in Subsection~\ref{sec:Glimsoft_full}.
	
	\begin{proposition}\label{stm:Glim-soft}
		For an infinitesimal sequence $\{\eps_k\}$, 
		consider $\cJ_k^0$ and $\cJ^0$ as in \eqref{eq:Jeps0} and \eqref{eq:J0}, respectively.
		Let also $P\in W^{1,q}(\Omega;\SL(3))$.
		\begin{enumerate}
			\item For every sequence
			$\{(v_k,P_k)\}\subset W^{1,2}_0(\Omega^0_k;\real^3) \times W^{1,q}(\Omega;\SL(3))$ such that
			$\{\eps_k \nv_k\}$ is $2$-equiintegrable and that
			$P_k \to P$ uniformly, we have
			\begin{equation*}
				\cJ^0(P) \leq \liminf_{k\to+\infty} \cJ_k^0(v_k,P_k).
			\end{equation*}
			\item There exists a bounded sequence $\{v_k\} \subset L^2(\Omega;\real^3)$, with $\{v_k\} \subset 	W^{1,2}_0(\Omega^0_{k};\real^3)$ for each $k$,
			such that
			\begin{equation*}
				\limsup_{k\to+\infty} \cJ_{k}^0(v_k,P_k) \leq \cJ^0(P),
			\end{equation*}
			provided $P_k \to P$ uniformly.
		\end{enumerate}
	\end{proposition}

	In the statement above,
	the space $W^{1,2}_0(\Omega^0_\eps;\real^3)$ is regarded for each $\eps$ as a subset of $W^{1,2}(\Omega;\real^3)$
	by extending its elements to $0$ on $\Omega^1_\eps$.
	The reason why we are only interested in functions with null traces roots in the splitting procedure,
	cf.~\eqref{eq:nulltrace} in Proposition~\ref{stm:splitting}.
	
	\begin{remark}
	Let $\Omega\subset \real^3$ be bounded Lipschitz domain and,
	for $p>1$, let us consider the local integral functionals on $W^{1,p}(\Omega;\real^3)$
	\begin{equation*}
		v \mapsto \io W_k(\nabla v) \dd x.
	\end{equation*}
	If the energy densities $\{W_k\}$ satisfy standard $p$-growth conditions,
	as a consequence of Rellich-Kondrachov theorem,
	the $\Gamma$-limits
	with respect to the strong $L^p$-convergence and
	with respect to the weak $W^{1,p}$-convergence
	coincide (if they exist).
	
	For the sequence of functionals
	\begin{equation}\label{eq:degeneratefunct}
		v \mapsto \io W_k(\eps_k \nabla v) \dd x,
	\end{equation}
	again under standard growth conditions for $\{W_k\},$   the analysis is more delicate.
	The natural bound that follows from the $p$-coercivity
	is $\| \eps_k \nabla v_k \|_{L^p} \leq C$,
	and it suggests the use of weak two-scale convergence
	(see Subsection~\ref{sec:twoscale}).
	However, this estimate alone is not enough to deduce convergence of the sequence $\{v_k\}$:
	a further control on the $\eps$-difference quotients is required
	to guarantee that a two-scale variant of Rellich-Kondrachov theorem holds
	(see \cite[Theorem~4.4]{Vis2}).
	
	In other words,   in our degenerate setting,
	compactness of sequences of gradients,
	say $\{\eps_k \nabla v_k\}$,
	does not bring compactness of $\{v_k\}$.
	This explains why we need to exploit the specific geometry of the perforated medium to recover the bound on $\|v_k\|_{L^2}$,
 see the proof of item (2) in Proposition~\ref{stm:Glim-soft}.
	
	We note incidentally that,
	by means of Lemma~\ref{lemma:2-scale}(4) below,
	it can be shown that the $\Gamma$-limit of the functionals \eqref{eq:degeneratefunct}
	with respect to the strong two-scale convergence in $L^p$ of $\{v_k\}$
	is the same as the one computed
	by combining the latter convergence and the weak two-scale convergence of $\{\eps_k \nabla v_k\}$.
	Those are not suitable choices for our goals, though,
	because, as we commented above,
	they do not match the natural compactness of the problem.
	This explains why in \cite{CC},
	where strong two-scale convergence is considered,
	the asymptotic behavior of minimum problems is not immediately determined
	by the $\Gamma$-convergence (see \cite[Section~10]{CC}).
	We also refer to Section~\ref{sec:conclusions} for a comparison
	between our findings and the ones in \cite{CC}.
	\end{remark}
	
\section{Preliminaries}\label{sec:pre}
We gather in this section the technical tools
to be employed in the sequel.

\subsection{A decomposition lemma}
	In our analysis of heterogeneous media it will be often desirable
	to disregard the energy contributions arising from the region close to $\partial \Omega$,
	for the composite fails to be periodic there
	(recall definitions \eqref{eq:Omega0}--\eqref{eq:Omega1}).
	To this aim,
	it is natural to resort to $p$-equiintegrability arguments,
	because such boundary strip has small measure.
	We recall that
	a family $\mathcal{C}\subset L^p(\Omega;\real^3)$ is said to be $p$-equiintegrable
	if for all $\delta>0$ there exists $m>0$ such that
	\[
	\sup_{u\in\mathcal{C}} \int_{E} \lvert u \rvert^p \dd x < \delta
	\quad\text{whenever } E \subset \Omega \text{ satisfies } \mathcal{L}^3(E)<m.
	\]
	
	The ensuing lemma grants that
	for any bounded sequence in $L^p$
	we can always find another one which is $p$-equiintegrable 
	and ``does not differ too much'' from the given one.
	
	\begin{lemma}[Theorem~2.20 in \cite{BaFo}; see also Lemma~1.2 in \cite{FMP}]
		\label{stm:dec-lemma}
		Let $\Omega$ be as in Section~\ref{sec:math}.
		For any sequence $\{v_k\} \subset W^{1,2}(\Omega;\real^3)$
		such that $v_k \rightharpoonup v$ weakly in $W^{1,2}(\Omega;\real^3)$
		there exist a subsequence $\{k_j\}$ and a sequence $\{u_j\} \subset W^{1,2}(\Omega;\real^3)$
		satisfying the following:
		\begin{enumerate}
			\item $u_j \rightharpoonup v$ weakly in $W^{1,2}(\Omega;\real^3)$;
			\item $u_j=v$ in a neighborhood of $\partial \Omega$;
			\item $\{\nabla u_j\}$ is $2$-equiintegrable;
			\item $\lim_{j\to+\infty} 
			\mathcal{L}^3 ( \{ x\in \Omega : v_{k_j} (x) \neq u_j(x) \} ) = 0.$
		\end{enumerate}
	\end{lemma}
	
	Property (4) yields
	$\lim_{j\to+\infty} 
	\mathcal{L}^3 (\{\nv_{k_j} \neq \nabla u_j \} )=0 $,
	because by standard properties of Sobolev functions
	(see, e.g., \cite[Lemma~7.7]{GiTr}) the inclusion
	$\{ v_{k_j} \neq u_j \} \supseteq \{\nv_{k_j} \neq \nabla u_j \}$ holds true.
	
\subsection{A couple of tools to deal with periodic heterogeneous media}
The periodic geometry of the composite calls for an extension result
for Sobolev maps on perforated domains.
Since the perforations of the matrix are well detached from the boundary,
by applying \cite[Lemma~B.7]{BrDFr} the following can be proved:
	
	\begin{lemma}[Lemma~8 in \cite{CC}]\label{stm:extension}
		Let $\Omega$ be open and bounded,
		and let $\Omega_\eps^1$ be as in Section~\ref{sec:intro}.
		There exists a linear and continuous extension operator
		\[
		\mathsf{T}_\eps \colon W^{1,2}(\Omega_\eps^1;\real^3)\to W^{1,2}(\Omega;\real^3)
		\]
		such that for all $y\in W^{1,2}(\Omega_\eps^1;\real^3)$
		\begin{gather*}
			\mathsf{T}_\eps y = y \quad \text{a.\,e. in } \Omega_\eps^1, \\
			\lVert \mathsf{T}_\eps y \rVert_{L^2(\Omega;\real^3)} \leq c\, \lVert y \rVert_{L^2(\Omega_\eps^1;\real^3)}, \\
			\lVert \nabla (\mathsf{T}_\eps y) \rVert_{L^2(\Omega;\matr)} \leq c\, \lVert \nabla y \rVert_{L^2(\Omega_\eps^1;\matr)},
		\end{gather*}
		where $c$ is independent of $\eps$ and $\Omega$.
	\end{lemma}
	
	\begin{remark}\label{stm:ext-omega}
		Even though the lemma above is a classical result,  
		it is worth clarifying the way we employ it.
		
		In the sequel, we always work with sequences
		which are already defined on the whole $\Omega$. 
		When we apply Lemma~\ref{stm:extension} to such a sequence,
		say $\{y_\eps\}\subset W^{1,2}(\Omega;\real^3)$,
		it is tacitly understood that
		the functions that are extended are the restrictions $y_\eps \llcorner \Omega_\eps^1$.
		So, in a sense, the process modifies $y_\eps$ on the region occupied by the soft inclusions
		rather than extending it.
		Note that the modification is a true one,
		because $\mathsf{T}_{\eps}$ cannot be the identity.
		The two crucial points for our analysis are that
		\begin{enumerate}
			\item if $\{y_\eps \llcorner \Omega_\eps^1\}$ and  $\{\ny_\eps \llcorner \Omega_\eps^1\}$ are bounded in $L^2$,
			then $\{\mathsf{T}_\eps y_\eps\}$ is bounded in $W^{1,2}(\Omega;\real^3)$;
			\item if $\{y_\eps\}$ is bounded in $L^2(\Omega;\real^3)$
			and $\{\ny_\eps\}$ is a $2$-equiintegrable sequence,
			then $\{\nabla (\mathsf{T}_\eps y_\eps)\}$ is $2$-equiintegrable as well.
		\end{enumerate}
		The second property follows from the construction of $\mathsf{T}_\eps$,
		which is modeled on the proof of \cite[Lemma~B.8]{BrDFr} by patching together
		the extensions from $W^{1,2}(Q^1;\real^3)$ to $W^{1,2}(Q;\real^3)$ given by \cite[Lemma~B.7]{BrDFr}
		via partitions of unity
		(this is also the reason why the constant $c$ above depends only on $Q^1$).
		The extensions in \cite[Lemma~B.7]{BrDFr} preserve equiintegrability,
		because they rely on the classical reflection procedure.	
	\end{remark}
	
	The first application of the extension lemma is the following Poincar\'e inequality on periodic heterogeneous media
	(cf.~formula~(4.5) in \cite{allaire} where, however, the proof is not provided).

	\begin{proposition}\label{stm:poincare}
		Let $\Omega$, $\Omega_\eps^0$ and $\Omega_\eps^1$ be as in Section~\ref{sec:intro}. There exists a constant $c$ independent of $\eps$ and
		such that for every $y \in W^{1,2}_0(\Omega;\real^3)$
		\begin{equation*}
			\| y \|_{L^2(\Omega;\real^3)}
			\leq c \left(
			\eps \| \nabla y \|_{L^2(\Omega^0_\eps;\matr)} +\| \nabla y \|_{L^2(\Omega^1_\eps;\matr)}
			\right).
		\end{equation*}
	\end{proposition}
	\begin{proof}
		For $\eps$ fixed, we use the extension operator $\mathsf{T}_\eps$ from Lemma~\ref{stm:extension} to obtain
		\begin{equation}\label{eq:poinc0}
			\begin{aligned}
				\| y \|_{L^2(\Omega)}
				& \leq \| y - \mathsf{T}_\eps y \|_{L^2(\Omega)} + \| \mathsf{T}_\eps y \|_{L^2(\Omega)} \\
				& = \| y - \mathsf{T}_\eps y \|_{L^2(\Omega^0_\eps)} + \| \mathsf{T}_\eps y \|_{L^2(\Omega)}. 
			\end{aligned}
		\end{equation}
		Observe that $\mathsf{T}_\eps y \in W^{1,2}_0(\Omega;\real^3)$, as $\mathsf{T}_\eps y = y$ a.\,e. in $\Omega_\eps^1$
		and there exists a tubular neighborhood $O$ of $\partial\Omega$ such that $\Omega^1_\eps \cap O \equiv \Omega \cap O$.	Then, by the standard Poincar\'e's inequality,
		\begin{equation}\label{eq:poinc1}
			\| \mathsf{T}_\eps y \|_{L^2(\Omega)}
			\leq c\| \nabla (\mathsf{T}_\eps y) \|_{L^2(\Omega)} 
			\leq c\| \nabla y \|_{L^2(\Omega^1_\eps)}.
		\end{equation}
		Observe also that $y - \mathsf{T}_\eps y \in W^{1,2}_0(\Omega^0_\eps;\real^3)$. In view of the periodic structure of $\Omega_\eps^0$ and of Poincar\'e inequality on each cube, we infer 
		\begin{align*}
			\|y - \mathsf{T}_\eps y\|_{L^2(\Omega^0_\eps)}^2 &= \sum_{t \in T_\eps} \|y - \mathsf{T}_\eps y\|_{L^2(\eps(t + Q_0))}^2 \\
			&= \sum_{t \in T_\eps} \eps^3 \int_{Q_0} |y(\eps(t + z)) - \mathsf{T}_\eps y(\eps(t + z))|^2 \dd z \\
			&\le c \sum_{t \in T_\eps} \eps^5 \int_{Q_0} |\nabla(y - \mathsf{T}_\eps y)(\eps(t + z))|^2 \dd z \\
			&= c \eps^2 \|\nabla(y - \mathsf{T}_\eps y)\|_{L^2(\Omega^0_\eps)}^2,
		\end{align*}
		where $c$ depends only on $Q_0$. By applying again Lemma~\ref{stm:extension} we find
		$$\|y - \mathsf{T}_\eps y\|_{L^2(\Omega^0_\eps)} \le c \eps \left(\|\nabla y\|_{L^2(\Omega^0_\eps)} + \|\nabla y\|_{L^2(\Omega^1_\eps)} \right).$$
		This, together with \eqref{eq:poinc0} and \eqref{eq:poinc1}, yields the result.
	\end{proof}
	
	\subsection{Two-scale convergence and the unfolding method}\label{sec:twoscale} 
	From a mathematical perspective,
	the high-contrast structure of the functional $\cJ_\eps$
	results in the absence of uniform bounds in $L^2$
	for sequences with equibounded energy;
	indeed, only bounds on $\{\eps \nabla y_\eps P^{-1}_\eps\}$ are available.
	Such degenerate bounds are conveniently dealt with
	by means of two-scale convergence \cite{allaire,nguetseng},
	whose definition we recall next.
	Hereafter, the subscript $\rm per$ denotes spaces of $Q$-periodic functions, e.g.
	\[
		W^{1,2}_{\rm per}(\real^3) \coloneqq \{
			u\in W^{1,2}_{\rm loc}(\real^3)
			: u(x+t)=u(x) \text{ a.e.\ for all } t\in \inte^3
		\}.
	\]
	\begin{definition}
		\label{stm:twoscaleconv}
		Let $\{\eps_k\}\subset (0,+\infty)$ be infinitesimal.
		A sequence	$\{y_k\}\subset L^2(\Omega;\real^3)$ \emph{weakly two-scale converges in $L^2$}
		to a function $y\in L^2(\Omega;L^2_\per(\real^3;\real^3))$ (notation: $y_k \wts y$)
		if for every $v\in L^2(\Omega;C_\per(\real^3;\real^3))$
		\[
		\lim_{k\to+\infty} \int_{\Omega} y_k(x)\cdot v \bigg(x,\frac{x}{\eps_k}\bigg) \dd x = \int_\Omega \int_Q y(x,z) \cdot v(x,z) \dd z \dd x.
		\]
		A sequence $\{y_k\}\subset L^2(\Omega;\real^3)$ \emph{strongly two-scale converges in $L^2$}
		to $y \in L^2(\Omega;L^2_\per(\real^3;\real^3))$ (notation: $y_k \stackrel{2}{\to} y$)	
		if $y_k \wts y$ in $L^2$ and 
		$\|y_k\|_{L^2(\Omega;\real^3)}\to \|y\|_{L^2(\Omega\times Q;\real^3)}$.
	\end{definition}
	
	Recalling that for $i=0,1$
	$\chi^i_k(x)=1$ if $x\in \Omega^i_k$ and $\chi^i_k(x)=0$ otherwise,
	an example of strong two-scale convergence is provided
	by the sequences $\{\chi^0_k\}$ and $\{\chi^1_k\}$.
	Indeed,
	\begin{equation}\label{stm:conv-chik}
		\chi^i_k \stackrel{2}{\to} \chi^i
		\quad \text{strongly two-scale in } L^2,
	\end{equation}
	where $\chi^i(x,z) \coloneqq \chi_{Q^i}(z)$ for all $(x,z)\in \Omega \times Q$.
	
	We collect in the next lemma some basic properties of two-scale convergence
	which we will resort to in the following.
	Proofs and more details can be found in \cite{allaire,Vis,Vis2}.
	
	\begin{lemma}\label{lemma:2-scale}
		Let $\{\eps_k\}\subset (0,+\infty)$ be infinitesimal and
		consider $\{y_k\}\subset L^2(\Omega;\real^3)$.
		\begin{enumerate}
			\item If $\{y_k\}$ is weakly two-scale convergent,
			then it is bounded in $L^2(\Omega;\real^3)$; conversely,
			if $\{y_k\}$ is bounded in $L^2(\Omega;\real^3)$,
			then it admits a weakly two-scale convergent subsequence.
			\item If $y_k \wts y$ weakly two-scale in $L^2$,
			then $y_k \rightharpoonup \int_Q y(\,\cdot\,,z)\dd z$ weakly in $L^2(\Omega;\real^3)$.
			\item If $y_k \wts y$ weakly two-scale in $L^2$
			and if $\{\chi_k\}\subset L^\infty(\Omega)$ is a bounded sequence that converges to $\chi \in L^\infty(\Omega)$ in measure,
			then $\chi_k y_k \wts \chi y$ weakly two-scale in $L^2$.
			\item Suppose that $\{y_k\}\subset W^{1,2}(\Omega;\real^3)$ and that
			$\{y_k\}$ and $\{\eps_k \nabla y_k\}$ are bounded in $L^2$.
			Then, there exists $y \in L^2(\Omega;W^{1,2}_\per(\real^3;\real^3))$
			such that, up to subsequences,
			$y_k \wts y$ and $\eps_k \nabla y_k \wts \nabla_z y$ weakly two-scale in $L^2$. 
		\end{enumerate}
	\end{lemma}
	
	Two-scale convergence in $L^2$ can be related to $L^2$ convergence by means of \emph{unfolding operator},
	which, for $\eps>0$, is the map 
	$\mathsf{S}_\eps\colon L^2(\Omega)\to L^2(\real^3 \times Q;\real^3)$
	defined as
	\begin{equation}\label{eq:unfolding}
		\mathsf{S}_\eps y(x,z) \coloneqq  \hat y \Big( \eps \left\lfloor \frac{x}{\eps}\right\rfloor + \eps z \Big),
	\end{equation}
	where $\hat y$ denotes the extension of $y$ by $0$ outside $\Omega$ and $\lfloor \,\cdot\, \rfloor$ is the floor function.
	
	\begin{lemma}\label{stm:unfolding}
		If $\{y_\eps\}\subset L^2(\Omega;\real^3)$ is bounded,
		the following hold:
		\begin{enumerate}
			\item $y_\eps \wts y$ weakly two-scale in $L^2$ if and only if
			$\mathsf S_\eps y_\eps \rightharpoonup y$ weakly in $L^2(\real^3\times Q;\real^3)$;
			\item $y_\eps \stackrel{2}{\to} y$ strongly two-scale in $L^2$ if and only if
			$\mathsf S_\eps y_\eps \to y$ strongly in $L^2(\real^3\times Q;\real^3)$.
		\end{enumerate}
		In addition, if $\{y_\eps\}$ is $2$-equiintegrable,
		the family of unfoldings $\{\mathsf{S}_\eps y_\eps\}$ is also $2$-equiintegrable on $\real^3 \times Q$.
		Lastly, if $y\in W^{1,2}(\Omega;\real^3)$, then
		$$
			\mathsf S_\eps(\eps\nabla y)(x,z) = \nabla_z (\mathsf S_\eps y )(x,z) \ \mbox{ a.e.\ in }\real^3 \times Q.
		$$
	\end{lemma}
	For a proof of Lemma~\ref{stm:unfolding} and for further reading on the unfolding operator
	we refer to \cite{Vis,Vis2,CDG02,CDG08}.
	
\subsection{Homogenization of connected media in finite plasticity}
We present a variant of \cite[Theorem~2.2]{DGP1}
that is instrumental in dealing with the analysis of the stiff matrix.
Its proof is an adaptation of the one in \cite{DGP1},
the most substantial difference being the use
of \cite[Theorem~19.1]{BrDFr} instead of \cite[Theorem~14.5]{BrDFr}.

Recalling that we chose $q>3$, 
we work in the space $W^{1,2}(\Omega;\real^3)\times W^{1,q}(\Omega;K)$
endowed with the topology $\tau$ characterized by 
		\begin{equation}\label{eq:tau}
			(y_k,P_k) \stackrel{\tau}{\to} (y,P)
			\quad\text{if and only if}
			\quad
			\begin{cases}
				y_k \to y &\text{strongly in } L^2(\Omega;\real^3), \\[1 mm]
				P_k \to P &\text{uniformly}.
			\end{cases}
	\end{equation}

\begin{theorem}
	\label{stm:homo-fin-plast2}
	Let $E$ be an open and connected set that is $Q$-periodic and that has Lipschitz boundary. For every $(y,P) \in W^{1,2}(\Omega;\real^3) \times W^{1,q}(\Omega;K)$,
	let
	$$
	\widetilde{W}(x,F) \coloneqq \chi_E(x)W^1(F), \qquad \widetilde{H}(x,P) \coloneqq \chi_E(x)H(P),
	$$
	and define
	\begin{equation}\label{tildeFeps}
		\cF_\eps(y,P) \coloneqq \io \widetilde{W}\left(\frac{x}{\eps},\ny(x) P^{-1}(x)\right) \dd x + \io \widetilde{H}\left(\frac{x}{\eps},P(x)\right) \dd x + \io |\nabla P(x)|^q \dd x,
	\end{equation}
	which we extend by setting 
	$$
	\cF_\eps(y,P)=+\infty
	\quad\text{on }
	\big[L^2(\Omega;\real^3)\times L^q(\Omega;\SL(3))\big]
	\setminus \big[W^{1,2}(\Omega;\real^3)\times W^{1,q}(\Omega;K)\big].
	$$
	If $W^1$ and $H$ satisfy \ref{E-growth}--\ref{E-lip}
	and \ref{H1}--\ref{H2}, respectively,
	then for all $(y,P) \in  L^2(\Omega;\real^3) \times L^q(\Omega;\SL(3))$
	the $\Gamma$-limit
	$$\cF(y,P) \coloneqq \Gamma(\tau)\mbox{-}\lim_{\eps \to 0} \cF_\eps(y,P)$$
	exists and we have that
	$$\cF(y,P) = \left\lbrace\begin{aligned}
		\displaystyle\io \Big(\widetilde W_\homo(\ny(x),P(x)) &+ \widetilde H_\homo(P(x)) + |\nP(x)|^q\Big) \dd x \\
		& \text{if } (y,P) \in W^{1,2}(\Omega;\real^3)\times W^{1,q}(\Omega;K), \\[3pt]
		+\infty \qquad\qquad\qquad\qquad &\text{otherwise in } L^2(\Omega;\real^3)\times L^q(\Omega;\SL(3)),
	\end{aligned}\right.$$
	where $\widetilde W_\homo \colon \matr \times K \to [0,+\infty)$ and $\widetilde H_\homo \colon K \to [0,+\infty)$ are defined as
	\begin{gather*}
		\widetilde W_\homo(F,G) \coloneqq \lim_{\lambda  \to +\infty} \frac{1}{\lambda ^3} \inf\left\lbrace \int_{(0,\lambda)^3} \widetilde{W}\big(x,(F+\ny(x))G^{-1}\big) \dd x : y \in W^{1,2}_0((0,\lambda )^3;\real^3) \right\rbrace, \\
		\widetilde H_\homo(F) \coloneqq \int_Q \widetilde{H}(z,F) \dd z.
	\end{gather*}
\end{theorem}

We observe that the theorem above is similar in spirit to homogenization results for perforated domains. The case at stake is however different, in that later we will deal with functions defined on the nonperforated domain $\Omega$. This simplifies the analysis, because it spares us the need to extend Sobolev maps valued at $\SL(3)$.

Thanks to Lemma~\ref{stm:dec-lemma},
we are able to refine the choice of recovery sequences for $\cF$.
This will come in handy in the proof of Corollary~\ref{stm:conv-min}.

\begin{coro}\label{stm:ref-recovery}
	Under the assumptions of Theorem~\ref{stm:homo-fin-plast2},
	for any $(y,P)\in W^{1,2}(\Omega;\real^3)\times W^{1,q}(\Omega;K)$
	there exists a recovery sequence $(y_k,P_k)$ for $\cF(y,P)$ satisfying the following:
	\begin{enumerate}
		\item $y_k \rightharpoonup y$ weakly in $W^{1,2}(\Omega;\real^3)$;
		\item $y_k=y$ in a neighborhood of $\partial \Omega$;
		\item $\{\ny_k\}$ is $2$-equiintegrable.
	\end{enumerate}
\end{coro}
\begin{proof}
	Let $\{(w_k,P_k)\}$ be a recovery sequence for $\cF(y,P)$
	as provided by Theorem~\ref{stm:homo-fin-plast2}.
	We apply Lemma~\ref{stm:dec-lemma} to $\{w_k\}$.
	We deduce the existence of sequences
	$\{k_j\}$ and $\{u_j\} \subset W^{1,2}(\Omega;\real^3)$ such that the sequence defined by
	$$
	y_k \coloneqq
	\begin{cases}
		u_j &\text{if }k=k_j \text{ for some } j\in \nat,\\
		y						& \text{otherwise}
	\end{cases}
	$$
	satisfies properties (1)--(3) and $(y_k,P_k) \stackrel{\tau}{\to} (y,P)$.
	Moreover
	\[
	\lim_{j\to+\infty} \mathcal{L}^3(N_j) = 0,
	\]
	where $N_j \coloneqq \{ x\in \Omega : w_{k_j} (x) \neq u_j(x) \}$.
	
	We are left to prove that
	$\{(y_k,P_k)\}$ satisfies the upper limit inequality.
	Loosely speaking, this is a consequence of the fact that 
	passing to a $2$-equiintegrable sequence ``does not increase the energy''.
	Upon passing to a subsequence, which we do not relabel, we can assume that $\{\cF_{k}(y_k, P_k)\}$ is convergent.
	We first focus on the elastic and hardening parts
	of the energy $\cF_{k_j}$. It holds that
	\begin{align*}
		&\int_{\Omega}
		\left[ W\left(\frac{x}{\eps_{k_j}}, \nabla w_{k_j} P^{-1}_{k_j}\right)
		+  H \left( \frac{x}{\eps_{k_j}},P_{k_j} \right)
		\right] \!\dd x \\
		&\quad = \int_{N_{j}} \left[
		W\left(\frac{x}{\eps_{k_j}}, \nabla w_{k_j} P^{-1}_{k_j}\right)
		+ H \left( \frac{x}{\eps_{k_j}},P_{k_j} \right)
		\right] \!\dd x \\
		& \quad\quad
		+ \int_{\Omega\setminus N_j}  \left[
		W\left(\frac{x}{\eps_{k_j}}, \nabla u_{j} P^{-1}_{k_j}\right) +
		H \left( \frac{x}{\eps_{k_j}},P_{k_j} \right)
		\right] \!\dd x \\
		&\quad\ge \int_{\Omega\setminus N_{j}}  \left[
		W\left(\frac{x}{\eps_{k_j}}, \nabla u_{j} P^{-1}_{k_j}\right) +
		H \left( \frac{x}{\eps_{k_j}},P_{k_j} \right)
		\right] \!\dd x,
	\end{align*}
	so that
	\begin{equation*}
		\begin{aligned}
			&\limsup_{j \to +\infty}
			\int_{\Omega}	\left[ W\left(\frac{x}{\eps_{k_j}}, \nabla w_{k_j} P^{-1}_{k_j}\right)
			+  H \left( \frac{x}{\eps_{k_j}},P_{k_j} \right)
			\right] \!\dd x \\
			&\quad\ge \limsup_{j \to +\infty} \int_{\Omega\setminus N_{j}}  \left[
			W\left(\frac{x}{\eps_{k_j}}, \nabla u_{j} P^{-1}_{k_j}\right) +
			H \left( \frac{x}{\eps_{k_j}},P_{k_j} \right)
			\right] \!\dd x \\
			&\quad=\limsup_{j \to +\infty}
			\int_{\Omega} \left[
			W\left(\frac{x}{\eps_{k_j}}, \nabla u_{j} P^{-1}_{k_j}\right) +
			H \left( \frac{x}{\eps_{k_j}},P_{k_j} \right)
			\right] \!\dd x,
		\end{aligned}
	\end{equation*}
	where the equality follows from the growth condition \ref{E-growth} and from the $2$-equiintegrability of $\{\nabla u_j\}$ (recall that $\sup_{k\in\nat} \|P^{-1}_k\|_\infty \le C$), together with the boundedness of $H$.
	Therefore, coming back to the full functional $\cF_{k_j}$,
	\begin{equation}\label{eq:limsup2}
		\begin{aligned}
			&\lim_{j \to +\infty} \cF_{k_j}(w_{k_j}, P_{k_j}) \\
			&\quad \ge \limsup_{j \to +\infty} \int_{\Omega}	\left[ 	
			W\left(\frac{x}{\eps_{k_j}}, \nabla w_{k_j} P^{-1}_{k_j}\right)
			+  H \left( \frac{x}{\eps_{k_j}},P_{k_j} \right)
			\right] \!\dd x
			+ \liminf_{j \to +\infty} \io |\nP_{k_j}|^q \dd x \\
			&\quad \ge \limsup_{j \to +\infty}  \int_{\Omega}	\left[ 	
			W\left(\frac{x}{\eps_{k_j}}, \nabla u_{j} P^{-1}_{k_j}\right)
			+  H \left( \frac{x}{\eps_{k_j}},P_{k_j} \right)
			\right] \!\dd x
			+ \liminf_{j \to +\infty} \io |\nP_{k_j}|^q \dd x \\
			&\quad\ge \lim_{j \to +\infty} \cF_{k_j}(u_{j}, P_{k_j}).
		\end{aligned}
	\end{equation}
	
	Recalling that $\{(w_k,P_k)\}$ is a recovery sequence
    and that we can assume $\{\cF_{k}(y_k, P_k)\}$ to be convergent,
	we find
	\begin{align*}
		\lim_{k \to +\infty} \cF_{k}(y_{k}, P_{k})
		= \lim_{j \to +\infty} \cF_{k_j}(u_{j}, P_{k_j})
		\le \lim_{j \to +\infty} \cF_{k_j}(w_{k_j}, P_{k_j})
		= \cF(y,P), 
	\end{align*}
	which in turn yields that $\{(y_k,P_k)\}$ is also a recovery sequence.
\end{proof}

\subsection{Finsler structure on $\SL(3)$}\label{sec:Finsler}
	In order to apply the results on homogenization
	of elastoplastic media in \cite{DGP1},
	we endow $\SL(3)$ with a Finsler structure.
	In doing so, we follow \cite{Mie02},
	whose approach is based on the notion of \emph{plastic dissipation}.
	Such line of thought links the geometry of $\SL(3)$ 
	to the physics of the system under consideration,
	and allows to conveniently include dissipation effects in the model,
	see Subsection~\ref{sec:diss}.
	
	We start from the observation that
	$\SL(3)$ is a smooth manifold
	with respect to the topology induced by the inclusion in $\matr$.
	For every $F \in \SL(3)$ the tangent space at $F$ is characterized as
	$$
	\rmT_F \SL(3)=F{\sf sl}(3) \coloneqq \{FM \in \matr : \tr M = 0\},
	$$
	and, in particular, $\rmT_I \SL(3)$ coincides with
	${\sf sl}(3)\coloneqq \{M \in \matr : \tr M = 0\}$.
	To the purpose of endowing $\SL(3)$ with a Finsler structure,
	we consider a function $\Delta_I \colon {\sf sl}(3)\to [0,+\infty)$
	on which we make the following assumptions (cf.~\cite[Section~1.1]{BaChSh} and \cite[Section~1]{Mie02}):
	\begin{enumerate}[label=\textbf{D\arabic*:},ref={D\arabic*}]
		\item\label{D0} It is $C^2$ on ${\sf sl}(3) \setminus \{0\}$;
		\item\label{D1} It is positively $1$-homogeneous:
		$\Delta_I(c M) = c\Delta_I(M)$ for all $c>0$ and $M \in {\sf sl}(3)$;
		\item\label{D3} The function $\Delta_I^2/2$ is strongly convex;
		\item\label{D2} It is $1$-coercive and has at most linear growth:
		there exist $0 < c_4 \le c_5$ such that for all $M \in {\sf sl}(3)$
		$$c_4|M| \le \Delta_I(M) \le c_5|M|.$$
	\end{enumerate}
	Note that we consider more restrictive
	regularity assumptions than the ones in \cite{Mie02}
	because we appeal to results of differential geometry,
	where smoothness is customarily required.
	The drawback of this choice is that
	in our analysis we cannot encompass some models, such as single crystal plasticity.
	However, on the positive side, our assumptions cover Von Mises plasticity,
	see \cite{HaMiMi,Mie02} and \cite[Example~3.6]{DGP1}.
	
	Let $ \rm T\SL(3)$ denote the tangent bundle to $\SL(3)$.
	We can ``translate'' $\Delta_I$
	to the tangent spaces other than ${\sf sl}(3)$
	by setting
	\begin{equation}\label{eq:Delta}
		\begin{array}{rccc}
			\Delta\colon & \rm T\SL(3) &\to& [0,+\infty)\\[1pt]
			& (F,M) &\mapsto& \Delta_I(F^{-1}M).
		\end{array}
	\end{equation}
	Then, it can be proved that $(\SL(3),\Delta)$ is a $C^2$ Finsler manifold.
	For an introduction to Finsler geometry
	we refer to the monograph \cite{BaChSh}.	
	
	Next, we introduce the family $\mathcal C(F_0,F_1)$ of piecewise $C^2$ curves
	$\Phi\colon [0,1] \to\SL(3)$
	such that $\Phi(0) = F_0$ and $\Phi(1) = F_1$.
	We set
	\begin{equation}\label{eq:diss-dist}
		D(F_0,F_1) \coloneqq
		\inf\left\lbrace
		\int_0^1 \Delta\big( \Phi(t),\dot{\Phi}(t) \big)\! \dd t 
		: \Phi \in \mathcal C(F_0,F_1)
		\right\rbrace,
	\end{equation}
	where $\dot{\Phi}$ is the velocity of the curve.
	The function $D$ provides a non-symmetric distance on $\SL(3)$:
	it is positive,
	attains $0$ if and only if it is evaluated on the diagonal of $\SL(3)\times\SL(3)$,
	and satisfies the triangular inequality;
	in general, however, $D(F_0,F_1)\neq D(F_1,F_0)$.
	
	Note that under assumptions \ref{D0}--\ref{D3} it follows that $\Delta$ is subadditive (see \cite[Theorem~1.2.2]{BaChSh}), hence convex. Therefore, by	an application of the direct method of the calculus of variations (cf.~\cite[Theorem~5.1]{Mie02})
	it can be proved that for every $F_0,F_1\in \SL(3)$ there exists a curve
	$\Phi \in C^{1,1}([0,1];\SL(3))$ such that $\Phi(0) = F_0$, $\Phi(1) = F_1$ and
	\begin{equation}
	D(F_0,F_1) = \int_0^1 \Delta\big(\Phi(t),\dot{\Phi}(t)\big) \dd t.
	\end{equation}
	We call such a $\Phi$ a \textit{shortest path} between $F_0$ and $F_1$.
%

A {\em geodesic} between $F_0$ and $F_1$, instead,
is a path that is a critical point of the length functional
under variations that do not alter the endpoints.
When for any couple of points in a given subset $S$ of a Finsler manifold
there is a unique shortest path contained in $S$ joining those points,
we say that $S$ is {\em geodesically convex}.

The existence of a compact set $K$ meeting the requirements in \ref{H1} is guaranteed by the following lemma, whose proof is the content of \cite[Remark~3.5]{DGP1}.
	 
	\begin{lemma}\label{stm:K-in-H2}
	Assume that a $C^2$ Finsler structure on $\SL(3)$ is assigned.
	Then, there exists a geodesically convex, compact neighborhood of $I$.
	\end{lemma}
	
	
\section{Compactness and splitting}\label{sec:comp+split}

We now turn to the analysis of the high-contrast energy in \eqref{eq:Jeps}. We investigate in this section the compactness properties of sequences with equibounded energy.
We will see that, as a consequence of the behavior of the hardening functional $H$, we can reduce the problem to the case of pure elasticity addressed by {\sc K. Cherdantsev \& M. Cherednichenko} \cite{CC},
and we adapt their approach.

\begin{lemma}[Compactness]\label{stm:cpt}
Let $\{\epsilon_k\}$ be an infinitesimal sequence.
We suppose that $\{(y_k,P_k)\}_{k\in\nat}\subset L^2(\Omega;\real^3) \times L^q(\Omega;\SL(3))$
satisfies
\begin{gather*}
	\|y_k\|_{L^2(\Omega;\real^3)} \le C, \qquad
	\cJ_k(y_k,P_k)\leq C
\end{gather*}
for some $C\geq 0$, uniformly in $k$.
Let us denote by $\tilde y_k$ the extension of $y_k$
in the sense of Remark~\ref{stm:ext-omega} above.
Then, there exist subsequences of
$\{\epsilon_k\}$, $\{y_k\}$, and $\{P_k\}$, which we do not relabel,
as well as
$y\in L^2(\Omega;W^{1,2}_\per(\real^3;\real^3))$,
$y^1\in W^{1,2}(\Omega;\real^3)$,
$v \in L^2(\Omega;W^{1,2}_0(Q^0;\real^3))$, and
$P\in W^{1,q}(\Omega;K)$ such that the following hold:
\begin{gather}
	\label{eq:struct-y}
	y(x,z) = y^1(x) + v(x,z) \quad\text{for a.\,e. } (x,z)\in \Omega\times Q, \\
	\label{c1}
	y_k \wts y \quad\text{weakly two-scale in } L^2, \\
	\label{c2}
	\eps_k \ny_k \wts \nabla_z v \quad\text{weakly two-scale in } L^2, \\
	\label{eq:tildeyk}
	\tilde y_k \rightharpoonup y^1 \quad\text{weakly in } W^{1,2}(\Omega;\real^3), \\
	\nonumber
	P_k \to P, \quad P_k^{-1}\to P^{-1}
	\quad \text{weakly in } W^{1,q}(\Omega;\SL(3)) \text{ and uniformly in } C(\bar{\Omega};\SL(3)), \\
	\label{eq:nyk-Fk}
	\nabla \tilde y_k P_k^{-1} \rightharpoonup \nabla y^1 P^{-1} \quad\text{weakly in } L^2(\Omega;\matr).
\end{gather}
\end{lemma}

\begin{proof}
From the definition of $\cJ_k$, for all $k\in\nat$
\begin{equation}\label{e3}
\left\|\nabla P_k \right\|_{L^q} \le C.
\end{equation}
Besides, for all $k$, hypothesis~\ref{E-growth-eps}, the definition of $H$ and the bound \eqref{Pbound} imply
\begin{gather}
\left\|\eps_k \chi^0_k \ny_k P^{-1}_k \right\|_{L^2} + \left\|\chi_k^1 \ny_k P_k^{-1} \right\|_{L^2} \le C, \label{e1} \\
\left\|P_k\right\|_{L^\infty} + \left\|P_k^{-1}\right\|_{L^\infty} \le C. \label{e2}
\end{gather}
Thanks to \eqref{stm:Pbound}, from the first estimate we deduce
\begin{equation}\label{e4}
\left\|\eps_k \chi_k^0\ny_k \right\|_{L^2}
+ \left\|\chi_k^1\ny_k \right\|_{L^2} \le C,
\end{equation}
which is precisely formula~(21) in \cite{CC}.
Thus, for what concerns the sequence of deformations,
the same bounds as the purely elastic case are retrieved.
While referring to \cite{CC} for details,
here we limit ourselves to sketch
how \eqref{e4} entails two-scale compactness.

The boundedness of $\{y_k\}$ in $L^2$ and Lemma~\ref{lemma:2-scale}(4) yield the existence of a function $y \in L^2(\Omega;W^{1,2}_\per(\real^3;\real^3))$ such that,
up to subsequences, \eqref{c1} holds and
\begin{gather}\label{eq:nablay}
\eps_k \ny_k \wts \nabla_z y \quad\text{weakly two-scale in } L^2.
\end{gather}
Thanks to \eqref{stm:conv-chik} and Lemma~\ref{lemma:2-scale}(3), we also infer that
\begin{gather*}
\chi_k^1 y_k \wts \chi^1 y,
\quad 
\eps_k \chi_k^1 \ny_k \wts \chi^1 \nabla_z y \quad\text{weakly two-scale in } L^2.
\end{gather*}
Moreover, there exist $y^1\in W^{1,2}(\Omega;\real^3)$ and $v \in L^2(\Omega;W^{1,2}_0(Q^0;\real^3))$
such that the decomposition \eqref{eq:struct-y} and the convergence \eqref{eq:tildeyk} hold.
By combining \eqref{eq:struct-y} and \eqref{eq:nablay}, \eqref{c2} follows.

We now turn to the sequence of plastic strains.
By \eqref{e3} and \eqref{e2}, we see that $\{P_k\}$ is bounded in $W^{1,q}(\Omega;K)$.
Since $q>3$, Morrey's embedding yields the uniform convergence
of (a subsequence of) $\{P_k\}$ to some $P\in W^{1,q}(\Omega;K)$ (note that the uniform convergence of $\{P_k\} \subset W^{1,q}(\Omega;K)$ yields that $P$ attains values in $K$ as well). Therefore, by definition of the inverse matrix
\begin{equation*}
	P_k^{-1} = \frac{\left({\rm cof}P_k\right)^T}{\det P_k} = \left({\rm cof}P_k\right)^T,
\end{equation*}
we also deduce that $P_k^{-1} \to P^{-1}$ uniformly.


Finally, we observe that,
thanks to \eqref{eq:tildeyk} and the uniform convergence of $\{P_k^{-1}\}$,
\eqref{eq:nyk-Fk} is also inferred.
\end{proof}

It is well-known that $\Gamma$-limits are not additive.
In our case, however, we are able to show that
the asymptotic behavior of the functionals $\cJ_\eps$ is given exactly
by the sum of the $\Gamma$-limits of the soft and of the stiff contributions.
Such splitting will enable us
to treat the $\Gamma$-limits of $\cJ_\eps^0$ and of $\cJ_\eps^1$ separately.
We premise a simple lemma,
which deals with the hardening part of the energy.
We recall that, for $i=0,1$,
$\chi_k^i$ is the characteristic function of  $\Omega^i_k$.

\begin{lemma}\label{stm:cont-hard}
Under assumptions \ref{H1}--\ref{H2},
for any sequence $\{P_k\} \subset W^{1,q}(\Omega;K)$ converging uniformly to $P\in W^{1,q}(\Omega;K)$	
it holds that
\begin{gather*}
\lim_{k\to+\infty} \io \chi_k^i(x) H\big( P_k(x) \big ) \dd x
= \mathcal{L}^3(Q^i) \io H\big( P(x) \big ) \dd x
\qquad\text{for $i=0,1$.}
\end{gather*}
\begin{proof}
Let us focus on the case $i=0$ first.
We set
$$
E^0 \coloneqq \bigcup_{t \in \mathbb{Z}^3} (t+Q^0)
= \real^3 \setminus \overline{E^1},
$$
By definition of $\Omega^0_k$ (see \eqref{eq:Omega0}), we have
$$
	\Omega \cap (\eps_k E^0 \setminus \Omega^0_k) \subset \{x\in \Omega : \mathrm{dist}(x,\partial \Omega)\leq \lambda\eps_k \}.
$$
Since $\{H(P_k)\}$ is uniformly bounded by \ref{H1} and \ref{H2},
we see that
\begin{align*}
	&\lim_{k\to+\infty} \io \chi_k^0(x) H\big( P_k(x) \big ) \dd x \\
	& \qquad = \lim_{k\to+\infty} \io \chi_{\eps_k E^0}(x) H\big( P_k(x) \big ) \dd x
	- \lim_{k\to+\infty} \io \big( \chi_{\eps_k E^0}(x) - \chi_k^0(x)\big) H\big( P_k(x) \big ) \dd x \\
	& \qquad = \lim_{k\to+\infty} \io \chi_{\eps_k E^0}(x) H\big( P_k(x) \big ) \dd x.
\end{align*}
Then, by the Lipschitz continuity of $H$ on its domain,
\begin{align*}
	\lim_{k\to+\infty} \io \chi_{\eps_k E^0}(x) H\big( P_k(x) \big) \dd x
	&=  \lim_{k\to+\infty} \io \chi_{\eps_k E^0}(x) H\big( P(x) \big) \dd x \\
	& =  \mathcal{L}^3(Q^0) \io H\big( P(x) \big ) \dd x.
\end{align*}

The case $i=1$ follows from the previous one
by the identities $\chi_k^1= \chi_\Omega - \chi_k^0$ and $\mathcal{L}^3(Q^1) = 1- \mathcal{L}^3(Q^0)$.
\end{proof}
\end{lemma}

The splitting process is explained by the ensuing proposition.

\begin{proposition}[Splitting]\label{stm:splitting}
Let $\{\epsilon_k\}$ be an infinitesimal sequence,
and let $\{(y_k,P_k)\}_{k\in\nat}\subset W^{1,2}(\Omega;\real^3) \times W^{1,q}(\Omega;\SL(3))$
be a sequence satisfying
\begin{gather*}
\|y_k\|_{L^2(\Omega;\real^3)} \le C, \quad
\cJ_k(y_k,P_k)\leq C
\end{gather*}
for some $C\geq 0$, uniformly in $k$.
Let $\tilde y_k$ be the extension of $y_k$
in the sense of Remark~\ref{stm:ext-omega},
and let $v \in L^2(\Omega;W^{1,2}_0(Q^0;\real^3))$ be as in Lemma~\ref{stm:cpt}.
Then, defining $v_k \coloneqq y_k - \tilde y_k$, the following hold:
\begin{gather}
	\label{eq:nulltrace}
	\{v_k\}\subset W^{1,2}_0(\Omega_k^0;\real^3), \\
	\nonumber
	\lVert v_k \rVert_{L^2(\Omega;\real^3)}\leq C, \\
	\label{eq:nvk-to-nzv}
	\eps_k \nv_k \wts \nabla_z v
	\quad\text{weakly two-scale in } L^2, \\
	\label{eq:split-liminf}
	\liminf_{k\to+\infty} \cJ_k^0(v_k,P_k) + \liminf_{k\to+\infty} \cJ_k^1(\tilde y_k,P_k)
	\leq \liminf_{k\to+\infty} \cJ_k(y_k,P_k), \\
	\label{eq:split-limsup}
	\limsup_{k\to+\infty} \cJ_k(y_k,P_k)
	\leq \limsup_{k\to+\infty} \cJ_k^0(v_k,P_k) + \limsup_{k\to+\infty} \cJ_k^1(\tilde y_k,P_k).
\end{gather}
Moreover, in \eqref{eq:split-liminf},
$\{v_k\}$ may be replaced
with another sequence $\{w_k\}\subset W^{1,2}_0(\Omega_k^0;\real^3)$
such that $\{\eps_k \nabla w_k\}$ is $2$-equiintegrable and
$\eps_k \nabla w_k \rightharpoonup 0$ weakly in $L^2(\Omega;\matr)$.
\end{proposition}
\begin{proof}
We first prove that
\eqref{eq:nvk-to-nzv} -- \eqref{eq:split-limsup} hold
for the sequence $\{v_k\}$.
Afterwards, we will show how to recover the equiintegrability for the sequence of gradients.

We split the functional $\cJ_k$ evaluated on $(y_k,P_k)$ as follows:
\begin{align}
\cJ_k(y_k,P_k) & = \cJ_k^0(y_k,P_k) + \cJ_k^1(y_k,P_k) \nonumber \\
& = \cJ_k^0(v_k,P_k) + \cJ_k^1(\tilde y_k,P_k) + \mathcal{R}_k, \label{eq:split}
\end{align}
where $\cJ_k^0$ and $\cJ_k^1$ are as in \eqref{eq:Jeps0} and \eqref{eq:Jeps1}, and
\begin{align*}
\mathcal{R}_k & \coloneqq
\cJ_k^0(y_k,P_k) - \cJ_k^0(v_k,P_k) \\
& = \io \chi_k^0 \left[
W^0_\eps\left(\eps_k \ny_k P_k^{-1} \right) - W^0_\eps\left(\eps_k \nv_k  P_k^{-1}\right)
\right]  \dd x.
\end{align*}
We next show that $\mathcal{R}_k$ is asymptotically negligible.

Hypothesis~\ref{E-lip-eps} yields
\begin{align}\label{eq:Rk}
\left|\mathcal{R}_k\right|
\le c_3 \io \chi_k^0
\left(1 + \left| \eps_k \ny_k P_k^{-1} \right| + \left| \eps_k \nv_k P_k^{-1} \right| \right)
\left| \eps_k \nabla\tilde{y}_k P_k^{-1} \right| \dd x.
\end{align}
Since $\{(y_k,P_k)\}$ is equibounded in energy,
the sequences $\{\eps_k  \chi_k^0 \ny_k P_k^{-1}\}$,
$\{\chi_k^1 \ny_k P_k^{-1}\}$, and $\{P_k^{-1}\}$ are bounded in suitable Lebesgue spaces
(see \eqref{e1} and \eqref{e2}).
By the properties of the extension operator $\mathsf{T}_\eps$ in Lemma~\ref{stm:extension},
we deduce that
\begin{align*}
\io \left| \nabla \tilde y_k P_k^{-1} \right|^2 \dd x
\leq c \io \left| \nabla \tilde y_k  \right|^2 \dd x
\leq c \io \left| \chi_k^1 \ny_k \right|^2 \dd x
\leq c \io \left| \chi_k^1 \ny_k P_k^{-1}\right|^2 \dd x
\leq C
\end{align*}
(recall estimate \eqref{stm:Pbound}).
So, thanks to \eqref{c2}, we deduce that 
$$
\eps_k \nv_k = \eps_k \ny_k - \eps_k \nabla \tilde y_k \wts \nabla_z v
\quad\text{weakly two-scale in } L^2,
$$
In particular, by Lemma~\ref{lemma:2-scale}(1),
$\{\eps_k \chi_k^0 \nv_k P_k^{-1}\}$ is bounded in $L^2(\Omega;\matr)$.
By applying H\"older's inequality to the right-hand side of \eqref{eq:Rk},
we then find $ \mathcal{R}_k = O(\eps_k)$.
Owing to \eqref{eq:split} we conclude that
\eqref{eq:split-liminf} and \eqref{eq:split-limsup} hold.

To complete the proof, we are only left to establish the existence of the sequence $\{w_k\}$.
Upon extraction of a subsequence,
which we do not relabel,
we may assume that
in \eqref{eq:split-liminf} the lower limit involving $\cJ^0_k$ is a limit.
From the equiboundedness of the energy,
by arguing as in the lines before \eqref{e4}, we get
\begin{equation}\label{eq:est-ny}
\|\eps_k \ny_k\|_{L^2} \le C, \qquad \|\chi_k^1 \ny_k\|_{L^2} \le C.
\end{equation}
Then, \eqref{c2} holds and,
by Lemma~\ref{lemma:2-scale}(2), we obtain
\begin{equation*}
		\eps_k \ny_k \rightharpoonup 0
		\quad\text{weakly in } L^2(\Omega;\matr).
\end{equation*}
Lemma~\ref{stm:dec-lemma}
applied to the sequence $\{\eps_k \ny_k\}$
yields two sequences, $\{k_j\}$ and $\{u_j\}\subset W^{1,2}(\Omega;\real^3)$,
such that $\{\eps_{k_j} \nabla u_j\}$ is $2$-equiintegrable,
\begin{gather}
	\label{eq:epsnu_conv}
	\eps_{k_j} \nabla u_j \rightharpoonup 0
	\quad \text{weakly in } L^2(\Omega;\matr),\\
	\nonumber
	\lim_{j\to+\infty} \mathcal{L}^3(N_j)=0,
	\quad\text{with }
	N_j\coloneqq \{ x\in \Omega : y_{k_j} (x) \neq u_j(x)\}.
\end{gather}
Besides, we have
\begin{gather}
	\eps_{k_j} \chi_{k_j}^1 \nabla u_j \to 0 \quad \text{ strongly in } L^2(\Omega;\matr). \label{eq:epschinu_conv}
\end{gather}
Indeed, it holds that
\begin{align*}
\|\eps_{k_j} \chi_{k_j}^1\nabla u_j\|_{L^2(\Omega)} &= \|\eps_{k_j} \chi_{k_j}^1\nabla u_j\|_{L^2(N_j)} + \|\eps_{k_j} \chi_{k_j}^1\nabla y_{k_j}\|_{L^2(\Omega\setminus N_j)} \\
&\le \|\eps_{k_j}\nabla u_j\|_{L^2(N_j)} + \eps_{k_j} \|\chi_{k_j}^1\nabla y_{k_j}\|_{L^2(\Omega)},
\end{align*}
and the conclusion follows
by the $2$-equiintegrability of $\{\eps_{k_j} \nabla u_j\}$ and
from \eqref{eq:est-ny}.

We now define $\tilde{u}_j \coloneqq \mathsf{T}_{k_j} u_j$, with $\mathsf{T}_{k_j}$ as in Lemma~\ref{stm:extension}.
From Remark~\ref{stm:ext-omega} it follows that
$\{\eps_{k_j} \nabla\tilde{u}_j\}$ is $2$-equiintegrable as well. 
Thus, the sequence
defined by
$$
	w_k \coloneqq
		\begin{cases}
		u_j - \tilde{u}_j &\text{if }k=k_j \text{ for some } j\in \nat,\\
		0						& \text{otherwise}
		\end{cases}
$$
has the properties that
$w_k\in W^{1,2}_0(\Omega_k^0;\real^3)$ and
$\{\eps_k \nabla w_k\}$ is $2$-equiintegrable.
Moreover,
$$\eps_k \nabla w_k \rightharpoonup 0 \quad \text{weakly in } L^2(\Omega;\matr).$$
To see this, we write
$$\eps_{k_j} \nabla w_{k_j} = \eps_{k_j} \nabla u_j - \eps_{k_j} \nabla\tilde{u}_j.$$
The first term converges to $0$ weakly in $L^2(\Omega;\matr)$, as stated in \eqref{eq:epsnu_conv}. Additionally, Lemma~\ref{stm:extension} entails
$$\|\eps_{k_j} \nabla\tilde{u}_j\|_{L^2} \le c \|\eps_{k_j} \chi^1_{k_j} \nabla u_j\|_{L^2},$$
and the weak convergence of $\{\eps_{k} \nabla w_k\}$ follows from \eqref{eq:epschinu_conv}.

We are now ready to prove the validity of \eqref{eq:split-liminf}
when $\{\eps_k \nv_k\}$ is replaced
by the $2$-equiintegrable sequence $\{\eps_k \nabla w_k\}$.
By the definition of the sequence at stake, we have
\begin{equation}\label{eq:v-minus-w}
	\eps_{k_j}(\nv_{k_j} - \nabla w_{k_j}) 
	= \eps_{k_j} (\ny_{k_j} - \nabla u_j ) - 
		\eps_{k_j} ( \nabla\tilde{y}_{k_j} - \nabla\tilde{u}_j )
	\quad \text{a.\,e. in } \Omega.
\end{equation}
Lemma~\ref{stm:extension} yields
\begin{align*}
\eps_{k_j} \|\nabla\tilde{y}_{k_j} - \nabla\tilde{u}_j\|_{L^2(\Omega)} &= \eps_{k_j} \|\nabla\big(\mathsf{T}_{k_j}(y_{k_j} - u_j)\big)\|_{L^2(\Omega)} \\
&\le c\eps_{k_j} \|\chi_{k_j}^1 \nabla(y_{k_j} - u_j)\|_{L^2(\Omega)} \\
&= c\eps_{k_j} \|\chi_{k_j}^1 (\ny_{k_j} - \nabla u_j)\|_{L^2(N_j)} \\
&\le c\left(\eps_{k_j} \|\chi_{k_j}^1\nabla y_{k_j}\|_{L^2(\Omega)} + \|\eps_{k_j}\nabla u_j\|_{L^2(N_j)}\right).
\end{align*}
Thus, \eqref{eq:est-ny} and the $2$-equiintegrability of $\{\eps_{k_j} \nabla u_j\}$ entail
\begin{equation}\label{eq:extens-conv}
\eps_{k_j} \left(\nabla\tilde{y}_{k_j} - \nabla\tilde{u}_j\right) \to 0 \quad \text{ strongly in } L^2(\Omega;\matr).
\end{equation}
Therefore, using \eqref{eq:v-minus-w} and the fact that the densities $W^0_{k_j}$ are bounded from below, we have
\begin{align*}
&	\io \chi_{k_j}^0(x)
			W^0_{k_j}\big(\eps_{k_j}\nv_{k_j}(x)P^{-1}_{k_j}(x)\big)
		\dd x \\
&	\quad =
		\int_{N_j}
			\chi_{k_j}^0(x)
			W^0_{k_j}\big(\eps_{k_j}\nv_{k_j}(x)P^{-1}_{k_j}(x)\big)
		\dd x \\
&	\quad\quad
		+ \int_{\Omega \setminus N_j}
				\chi_{k_j}^0(x)
				W^0_{k_j}\Big(
					\big(
					\eps_{k_j}\nabla w_{k_j}(x)
					- \eps_{k_j}(\nabla\tilde{y}_{k_j}(x)
					- \nabla\tilde{u}_j(x))
					\big)P^{-1}_{k_j}(x)
				\Big)
				\dd x \\
&	\quad\quad
		- \int_{\Omega\setminus N_j}
			\chi_{k_j}^0(x)
			W^0_{k_j}\big(\eps_{k_j}\nabla w_{k_j}(x)P^{-1}_{k_j}(x)\big)
		\dd x 
		+ \int_{\Omega\setminus N_j}
			\chi_{k_j}^0(x)
			W^0_{k_j}\big(\eps_{k_j}\nabla w_{k_j}(x)P^{-1}_{k_j}(x)\big)
		\dd x \\
&\quad\ge - c\left(\int_{\Omega\setminus N_j} |\eps_{k_j}(\nabla\tilde{y}_{k_j}(x) - \nabla\tilde{u}_j(x))|^2 \dd x\right)^{1/2} + \int_{\Omega\setminus N_j} \chi_{k_j}^0(x)
W^0_{k_j}\big(\eps_{k_j}\nabla w_{k_j}(x)P^{-1}_{k_j}(x)\big) \dd x,
\end{align*}
where the Lipschitz regularity \ref{E-lip-eps} and H\"older's inequality were employed
to derive the last bound
(recall that $\sup_{k\in\nat} \|P^{-1}_k\|_\infty \le C$).
We now take the limit in the inequality above.
According to Lemma~\ref{stm:cont-hard}, the hardening term has a limit.
Therefore, also the elastic contribution is convergent,
and it satisfies
$$
\lim_{k\to+\infty} \cJ^0_k(v_k,P_k)
= \lim_{j\to+\infty} \io \chi_{k_j}^0(x)
W^0_{k_j}\big(\eps_{k_j}\nv_{k_j}(x)P^{-1}_{k_j}(x)\big) \dd x
+ \mathcal{L}^3(Q^0) \io H\big( P(x) \big ) \dd x.
$$
The strong converge \eqref{eq:extens-conv} implies
\begin{align*}
&\lim_{j\to+\infty} \io \chi_{k_j}^0(x)
W^0_{k_j}\big(\eps_{k_j}\nv_{k_j}(x)P^{-1}_{k_j}(x)\big) \dd x \\
&\quad\ge \liminf_{j\to+\infty} \int_{\Omega\setminus N_j} \chi_{k_j}^0(x)
W^0_{k_j}\big(\eps_{k_j}\nabla w_{k_j}(x)P^{-1}_{k_j}(x)\big)  \dd x \\
&\quad = \liminf_{j\to+\infty} \io \chi_{k_j}^0(x)
W^0_{k_j}\big(\eps_{k_j}\nabla w_{k_j}(x)P^{-1}_{k_j}(x)\big) \dd x,
\end{align*}
where the equality follows from the growth condition \ref{E-growth-eps} and
from the equiintegrability of $\{\eps_{k_j}\nabla w_{k_j}\}$.
We thereby infer
$$
\liminf_{k\to+\infty} \cJ_k^0(w_k,P_k)
\leq \liminf_{j\to+\infty} \cJ_{k_j}^0(w_{k_j},P_{k_j})
\leq \lim_{k\to+\infty} \cJ_{k}^0(v_{k},P_{k}),
$$
and this concludes the proof.
\end{proof} 


\section{$\Gamma$-limit of the soft component}\label{sec:Glimsoft}
We devote this section to the study of the asymptotics
of the functional $\cJ^0_\eps$ in \eqref{eq:Jeps0},
which encodes the energy of the inclusions.
After some observations on the limiting functional $\cJ^0$ in \eqref{eq:J0},
in the second and third subsections we deal respectively
with the lower and with the upper limit inequality for the elastic part of the energy.
The other contributions will be taken into account in Subsection~\ref{sec:Glimsoft_full},
where we prove Proposition~\ref{stm:Glim-soft}.

\subsection{The limiting functional}
The definition of $\rmQ W^0$ in \eqref{eq:tildeW0},
which encodes the limiting elastic contribution of the soft inclusions,
may be regarded as a variant
of the well known Dacorogna's formula for the quasiconvex envelope
\cite[Theorem~6.9]{dacorogna}.
As such, the infimum in \eqref{eq:tildeW0} does not depend on $Q$,
and we may rewrite $\rmQ W^0$ as follows:
\begin{equation}\label{eq:Q'W-bis}
	\rmQ W^0(F,G) = \inf \left\{ \fint_{Q^0}
	W^0 \Big(\big( F + \nv(z) \big) G \Big) \dd z : v\in W^{1,2}_0(Q^0;\real^3)\right\}.
\end{equation}

Note that here quasiconvexification occurs just with respect to the first argument,
since a very strong convergence is considered for the second one (cf.~Proposition~\ref{stm:Gliminf-soft} below).
The fact that
different variables in a problem may call for different relaxation procedures
has been already observed.
As an example, we mention
the concept of cross-quasiconvexity introduced by {\sc Le~Dret \& Raoult} \cite{LeDRa}
to deal with dimension reduction problems in elasticity.
 
For the sake of completeness,
we explicitly mention some basic properties of $\rmQ W^0$.
Note that in the following lemma we use the symbol $W^0$ to denote a generic function rather than the specific one in \ref{E-conv}.

\begin{lemma}\label{stm:lip}
Let $W^0\colon \matr \to \real$, and assume there exist $0 < c_1 \le c_2$ such that for all $F\in\matr$
$$
	c_1 |F|^2\leq W^0(F) \le c_2\left(|F|^2+1\right).
$$
Let	$\rmQ W^0$ be as in \eqref{eq:tildeW0}.
\begin{enumerate}
\item For all $F,G\in\matr$
$$
	c_1 |FG|^2\leq \rmQ W^0(F,G) \le c_2\left(|FG|^2+1\right),
$$
and
for all $G\in\matr$ there exists $c\coloneqq c(G) > 0$ such that for all $F_1,F_2 \in \matr$
\begin{align*}
\left|\rmQ W^0(F_1,G)-\rmQ W^0(F_2,G)\right| \le c \left(1+|F_1|+|F_2|\right)|F_1-F_2|.
\end{align*}
\end{enumerate}
Suppose further that
there exists $c_3 > 0$ such that for all $F_1,F_2 \in \matr$
\begin{equation}\label{eq:W0-lip}
\left|W^0(F_1)- W^0(F_2)\right| \le c_3 \left(1+|F_1|+|F_2|\right)|F_1-F_2|.
\end{equation}
\begin{enumerate}[resume]
\item Then, $\rmQ W^0(F,\,\cdot\,)$ is continuous for all $F\in\matr$.
\item If $\{P_k\}\subset W^{1,q}(\Omega;\SL(3))$
converges weakly to $P\in W^{1,q}(\Omega;\SL(3))$,
then for any $V\in L^2(\Omega;\matr)$
\[
\lim_{k\to+\infty} \io \rmQ W^0\big( V(x), P^{-1}_k(x) \big) \dd x
= \io \rmQ W^0\big( V(x), P^{-1} (x) \big) \dd x.
\]
\end{enumerate}
\end{lemma}
\begin{proof}
The growth conditions on $\rmQ W^0$ are an immediate consequence
of the ones on $W^0$ and of the definition of $\rmQ W^0$.

For what concerns the $2$-Lipschitz property,
let us set $W^0_G(F)\coloneqq W^0(FG)$ for any fixed $G\in\matr$.
Then, $\rmQ W^0(\,\cdot\,,G)$ coincides with
the quasiconvex envelope of $W^0_G$.
By \cite[Remark~5.4(iii)]{dacorogna} it follows that $\rmQ W^0(\,\cdot\,,G)$ is separately convex,
and hence, in view of the growth assumptions on $W^0$,
the proof of item (1) is concluded by \cite[Proposition~2.32]{dacorogna}.

As for assertion (2), let $G_k \to G$ in $\matr$.
In view of \eqref{eq:W0-lip}, 
for every $\delta>0$ there exists $c_\delta>0$ such that
$$
\rmQ W^0(F,G_k)  - \rmQ W^0(F,G) \leq c_\delta | G_k - G |+\delta. 
$$
Similarly, for any $k\in\nat$ there exists $v_k\in W^{1,p}_0(Q;\matr)$ such that
\begin{multline*}
\rmQ W^0(F,G_k)  - \rmQ W^0(F,G)
\\ \geq - c_3 | G_k - G |\int_Q \big( 1+ |(F+\nv_k)G_k|+|(F+\nv_k)G|\big) |F+\nv_k|\dd x - \frac{1}{k}.
\end{multline*}
Thanks to the coercivity of the integrand, it follows that
$\{\nv_k\}$ is bounded in $L^2$, whence
$$
\rmQ W^0(F,G_k)  - \rmQ W^0(F,G) \geq - c\, | G_k - G |- \frac{1}{k} 
$$
for a constant $c$ independent of $k$.
The continuity of $\rmQ W^0(F,\,\cdot\,)$ is then proved
by letting first $k\to+\infty$ and then $\delta\to0$.

Finally, taking into account properties (1) and (2),
as well as the compact embedding of $W^{1,q}(\Omega)$ into $C(\bar\Omega)$,
we can employ the dominated convergence theorem
to obtain the continuity property in (3).
\end{proof}

We now exhibit an alternative expression for the soft limiting elastic energy,
which is to be exploited in the proof of Proposition~\ref{stm:Glimsup-soft}.

\begin{lemma}\label{stm:int-tildeW}
For every couple $(V,P)\in L^2(\Omega;\matr) \times W^{1,q}(\Omega;\SL(3))$
we have
\begin{equation}\label{eq:int-tildeW}
\begin{aligned}
&\io \rmQ W^0\big(V(x),P^{-1}(x)\big) \dd x
\\ &= \inf\left\{
\io \fint_{Q^0}
W^0 \Big(\big( V(x) + \nabla_z w(x,z) \big) P^{-1}(x) \Big) \dd z \dd x : w\in L^2(\Omega;W^{1,2}_0(Q^0;\real^3))
\right\}.
\end{aligned}
\end{equation}
\end{lemma}
The identity above rests on a measurable selection criterion
that we recall next.
\begin{lemma}[Lemma~3.10 in \cite{FoKr}]
\label{stm:FoKr}
Let $S$ be a multifunction
defined on the measurable space $X$
and taking values in the collection of subsets of the separable metric space $Y$.
If $S(x)$ is nonempty and open in $Y$ for every $x\in X$,
and if the set $\Set{ x \in X : y\in S(x)}$ is measurable for every $y\in Y$,
then $S$ admits a measurable selection,
that is,
there exists a measurable function
$s\colon X \to Y$ such that $s(x) \in S(x)$ for all $x\in X$.
\end{lemma}
The previous lemma is a variant of \cite[Theorem~III.6]{CaVa},
and we refer to that monograph for a comprehensive treatment of measurable selection principles.

\begin{proof}[Proof of Lemma~\ref{stm:int-tildeW}]
The argument follows the one proposed in \cite[Corollary~3.2]{FoKr}.

Let us fix $w\in L^2(\Omega;W^{1,2}_0(Q^0;\real^3))$,
so that, for almost every $x\in \Omega$,
$w(x,\,\cdot\,)\in W^{1,2}_0(Q^0;\real^3)$.
Hence, according to \eqref{eq:Q'W-bis}, we have
$$
\rmQ W^0\big(V(x),P^{-1}(x)\big)
\leq \fint_{Q^0} W^0 \Big(\big( V(x) + \nabla_z w(x,z) \big) P^{-1}(x) \Big) \dd z
\quad\text{ for a.\,e. }x\in \Omega,
$$
whence, after integration over $\Omega$, we deduce that
in \eqref{eq:int-tildeW} the left-hand side is smaller that the righ-hand one.

In order to establish the opposite inequality,
we first observe that, by \eqref{eq:Q'W-bis}, for every $x\in\Omega$ and every $\delta>0$
there exists $v_{x,\delta}\in W^{1,2}_0(Q^0;\real^3)$ such that
\begin{equation}\label{eq:vxdelta}
\fint_{Q^0} W^0 \Big(\big( V(x) + \nv_{x,\delta}(z) \big) P^{-1}(x) \Big) \dd z
- \rmQ W^0\big(V(x),P^{-1}(x)\big) < \delta.
\end{equation}
We introduce the multifunction $S$ defined for $x\in\Omega$ by
$$
S(x)\coloneqq \left\{ v \in W^{1,2}_0(Q^0;\real^3) : \eqref{eq:vxdelta} \text{ holds for } v_{x,\delta}=v \right\}.
$$
We show that it admits a measurable selection.
To this purpose observe that,
as a consequence of the growth assumptions on $W^0$
and the dominated convergence theorem,
$S(x)$ is a nonempty, open subset of $W^{1,2}_0(Q^0;\real^3)$
for every $x\in\Omega$.
Second, for every $v \in W^{1,2}_0(Q^0;\real^3) $
the set $\{ x \in \Omega : v\in S(x)\}$ is measurable,
because it is the sublevel set of a measurable function.

Thanks to Lemma~\ref{stm:FoKr},
for every $\delta>0$ we retrieve a measurable function $w_\delta\colon \Omega \to W^{1,2}_0(Q^0;\real^3)$
that  satisfies
\begin{equation*}
\io \fint_{Q^0}  W^0 \Big(\big( V(x) + \nabla_z w_\delta(x,z) \big) P^{-1}(x) \Big) \dd z \dd x
\leq \io \rmQ W^0\big(V(x),P^{-1}(x)\big) + O(\delta).
\end{equation*}
In particular, by the growth conditions on $W^0$, $w_\delta$ must belong to $L^2(\Omega;W^{1,2}_0(Q^0;\real^3))$.
Therefore, since $\delta$ is arbitrary, we conclude that
the left-hand side in \eqref{eq:int-tildeW} bounds from above the right-hand one.
\end{proof}


\subsection{Lower bound for the elastic energy}
The goal of this subsection is to prove the ensuing:
\begin{proposition}\label{stm:Gliminf-soft}
Let $\{W^0_k\}_{k}$ satisfy assumptions \ref{E-growth-eps}--\ref{E-conv}, and let $P \in  W^{1,q}(\Omega;\SL(3))$.
For every sequence
$\{(v_k,P_k)\}\subset W^{1,2}_0(\Omega^0_k;\real^3) \times W^{1,q}(\Omega;\SL(3))$
such that $\{\eps_k\nv_k\}$ is $2$-equiintegrable and $P_k \to P$ uniformly,
it holds that
\begin{equation}\label{eq:Gliminf-soft}
\mathcal{L}^3(Q^0) \io \rmQ W^0\big(0,P^{-1}(x)\big) \dd x \leq \liminf_{k\to+\infty} \io \chi_k^0(x)
W^0_k\big(\eps_k\nv_k(x)P^{-1}_k(x)\big) \dd x.
\end{equation}
\end{proposition}

At a first glance, it may look bizarre that
no convergence for the sequence  $\{\eps_k\nv_k\}$ is prescribed.
The statement becomes clearer
once we recall that
if $\mathcal Q f$ is the quasiconvex envelope of $f\colon \matr \to \real$, then
$$
\mathcal Q f(0) \leq \fint_\Omega f\big( \nabla v(x) \big) \dd x
$$
for any $v\in W^{1,\infty}_0(\Omega;\real^3)$.

In order to establish \eqref{eq:Gliminf-soft},
it is convenient to unfold the elastic energy.

\begin{lemma}\label{stm:unfold-en}
Let $\{W^0_k\}_{k}$ satisfy assumptions \ref{E-growth-eps}--\ref{E-conv}. For any $(v,P)\in W^{1,2}(\Omega;\real^3) \times W^{1,q}(\Omega;\SL(3))$ it holds that
\begin{equation}\label{eq:elastic}
\io \chi_k^0(x)
W^0_k\left(\eps_k\nv(x)P^{-1}(x)\right) \dd x
=
\sum_{t \in T_k}
\int_{\eps_k(t+Q)} \int_{Q^0} W^0_k \left(
\nabla_z \hat{v} (x,z) \hat{P}^{-1} (x,z)
\right) \dd z \dd x,
\end{equation}
where
$\hat{v} \coloneqq \mathsf{S}_k v$,
$\hat{P} \coloneqq \mathsf{S}_k P$
and $\mathsf{S}_k \coloneqq \mathsf{S}_{\varepsilon_k}$ is the unfolding operator
introduced in Lemma~\ref{stm:unfolding}.
%
\end{lemma}
\begin{proof}
According to the definition of $\Omega^0_k$ in \eqref{eq:Omega0},
the left-hand side of \eqref{eq:elastic} equals
$$
\eps_k^3 \sum_{t \in T_k}
\int_{Q^0} W^0_k \Big(\eps_k \nv \big( \eps_k(t+z) \big) P^{-1}\big( \eps_k(t+z) \big)\Big) \dd z.
$$
We use the unfolding operator
to rewrite this quantity as a double integral.
Recalling Lemma~\ref{stm:unfolding},
we first observe that 
for every $t \in T_k$ and $z \in Q^0$ we have the identities
\begin{gather*}
	\mathsf{S}_k(\eps_k\nv)(\eps_k t,z)
	= \eps_k\nv\big( \eps_k(t+z) \big) , 
	\quad
	\mathsf{S}_k P^{-1}(\eps_k t,z)
	= P^{-1}\big( \eps_k(t+z) \big) .
\end{gather*}
Then, we also have
\begin{gather*}
	\mathsf{S}_k(\eps_k\nv) = \nabla_z \big(\mathsf{S}_k v\big) = \nabla_z \hat{v},
	\quad
	\mathsf{S}_kP^{-1} = (\mathsf{S}_k P)^{-1} = \hat{P}^{-1}.
\end{gather*}
We obtain
\begin{align*}
\nonumber
& \io \chi^0_k(x) W^0_k \big( \eps_k\nv(x) P^{-1}(x) \big)\dd x \\
\nonumber
& \quad = \eps_k^3 \sum_{t \in T_k}
\int_{Q^0} W^0_k \big(
\mathsf{S}_k (\eps_k\nv) (\eps_k t,z)
\mathsf{S}_k (P^{-1}) (\eps_k t,z)
\big) \dd z \\
&\quad = \sum_{t \in T_k}
\int_{\eps_k(t+Q)} \int_{Q^0} W^0_k \left(
\nabla_z \hat{v} \left(\eps_k\left\lfloor\frac{x}{\eps_k}\right\rfloor,z \right)
\hat{P}^{-1} \left(\eps_k \left\lfloor\frac{x}{\eps_k}\right\rfloor,z\right)
\right) \dd z \dd x,
\end{align*}
because $\lfloor x / \eps_k \rfloor = t$ for all $x \in \eps_k(t+Q)$.
Since, in general, it holds that
$$\mathsf{S}_k u  \left( \eps_k \left\lfloor\frac{x}{\eps_k} \right\rfloor , z \right)
= u \left( \eps_k \left\lfloor\frac{x}{\eps_k} \right\rfloor + \eps_k z \right)
= \mathsf{S}_k u (x,z),$$
identity \eqref{eq:elastic} follows.
%
\end{proof}

A crucial ingredient in the proof of Proposition~\ref{stm:Gliminf-soft} is a sort of lower semicontinuity result for the elastic contribution to the energy.

\begin{lemma}\label{stm:lem20CC}
Let $\{W^0_k\}_{k}$ satisfy assumptions \ref{E-growth-eps}--\ref{E-conv}.
Let also $\{w_k\}\subset L^2(\Omega; W^{1,2}_0(Q_0;\real^3))$ be such that
$\{\nabla_z w_k\}$ is $2$-equiintegrable.
Then, for all $P \in W^{1,q}(\Omega;\SL(3))$,
$$
\mathcal{L}^3(Q^0) \io \rmQ W^0\big(0,P^{-1}(x)\big) \dd x
\le \liminf_{k\to+\infty} \int_{\Omega}\int_{Q^0} W^0_k \big( \nabla_z w_k(x,z ) P^{-1}_k(x) \big) \dd z \dd x,
$$
whenever $P_k \to P$ uniformly.
\end{lemma}

\begin{proof}
From \eqref{eq:Q'W-bis}
it follows that for all $k\in\nat$
\begin{align}\label{eq:step1}
\mathcal{L}^3(Q^0) \io \rmQ W^0\big(0,P^{-1}_k(x)\big) \dd x
\leq \io \int_{Q^0}
W^0 \big( \nabla_z w_k(x,z) P^{-1}_k(x) \big) \dd z \dd x.
\end{align}
%
%
%

Next, relying on the pointwise convergence of $\{W^0_k\}$ to $W^0$,
we adapt the argument in the proof of \cite[Theorem~5.14]{DalM}
to pass from $W^0$ to $W^0_k$ on the right-hand side
(see also \cite[Lemma~5.2]{DKP} for a similar result in the context of $\mathscr A$-quasiconvexity).	
Fix $\delta > 0$. If $\{\nabla_z w_k\}$ is $2$-equiintegrable, then so is $\{\nabla_z w_k P^{-1}_k\}$.
Therefore, since the $2$-growth assumptions on $\{W^0_k\}$ transfer to the pointwise limit $W^0$,
there exists $r > 0$ such that
\begin{equation}\label{eq:first-part}
\sup_{k \in\nat} \int_{\{(x,z) \in \Omega \times Q^0 : |\nabla_z w_k(x,z)P^{-1}_k(x)|>r\}} W^0\big( \nabla_z w_k(x,z) P^{-1}_k(x) \big) \dd z \dd x \le \delta.
\end{equation}
Owing to assumption \ref{E-lip-eps} and Remark~\ref{regW0}, we can find $\rho > 0$ such that
for every $F,G \in \matr$ contained in the open ball $B(0,\rho)$
\begin{equation}\label{eq:second-part}
\sup_{k\in\nat}|W^0_k(F) - W^0_k(G)| + |W^0(F) - W^0(G)| \le \delta.
\end{equation}
Let now $F_1,\dots,F_n \in B(0,r)$ be such that
\begin{equation}\label{eq:cover}
\overline{B(0,r)} \subset \bigcup_{i=1}^n B\left(F_i,\rho\right).
\end{equation}
Due to the pointwise convergence of $W^0_k$ to $W^0$, for any such $F_i$ there exist $\bar k_i \in \nat$ such that $|W^0_k(F_i) - W^0(F_i)| \le \delta$ if $k > \bar k_i$. Letting $\bar k \coloneqq \max\{\bar k_1,\dots,\bar k_n\}$, it follows that for any $i=1,\dots,n$
\begin{equation}\label{eq:unifconv}
|W^0_k(F_i) - W^0(F_i)| \le \delta \quad	\text{if } k> \bar k.
\end{equation}
By \eqref{eq:cover}, for every $G \in \overline{B(0,r)}$
there exists $i \in \{1,\dots,n\}$ such that $G \in B(F_i,\rho)$. For this particular $i$, the combination of the triangle inequality, \eqref{eq:second-part} and \eqref{eq:unifconv} yields
\begin{equation}\label{eq:lipschitz-est}
|W^0_k(G) - W^0(G)| \le |W^0_k(G) - W^0_k(F_i)| + |W^0_k(F_i) - W^0(F_i)| + |W^0(G) - W^0(F_i)| \le 3\delta,
\end{equation}
for every $G \in \overline{B(0,r)}$ and every $ k> \bar k$.

Thanks to Lemma~\ref{stm:lip}(3) and \eqref{eq:step1} we deduce
\begin{align*}
&\mathcal{L}^3(Q^0) \io \rmQ W^0\big(0,P^{-1}(x)\big) \dd x \\
&\qquad = \mathcal{L}^3(Q^0) \lim_{k\to+\infty}\io \rmQ W^0\big(0,P^{-1}_k(x)\big) \dd x \\
&\qquad\le \liminf_{k\to+\infty} \io \int_{Q^0}  W^0 \big( \nabla_z w_k(x,z) P^{-1}_k(x) \big) \dd z \dd x \\
&\qquad\le \liminf_{k\to+\infty} \int_{\{(x,z) \in \Omega \times Q^0 :|\nabla_z w_k(x,z)P^{-1}_k(x)|\leq r\}} W^0\big( \nabla_z w_k(x,z) P^{-1}_k(x) \big) \dd z \dd x + \delta \\
&\qquad\le \liminf_{k\to+\infty}\int_{\Omega}\int_{Q^0} W^0_k\big( \nabla_z w_k(x,z) P^{-1}_k(x) \big) \dd z \dd x + 3\delta\mathcal{L}^6(\Omega \times Q^0) + \delta,
\end{align*}
where the second inequality is due to \eqref{eq:first-part}, and the last one to \eqref{eq:lipschitz-est}.
The arbitrariness of $\delta > 0$ yields the conclusion.
\end{proof}

We are now ready to prove the lower bound for the elastic contribution of the soft part. 

\begin{proof}[Proof of Proposition~\ref{stm:Gliminf-soft}]
Let $\hat{v}_k \coloneqq \mathsf{S}_k v_k$.
In view of the $2$-equiintegrability of the sequence $\{\eps_k\nv_k\}$ and of Lemma~\ref{stm:unfolding},
$\{\nabla_z \hat{v}_k\}$ is $2$-equiintegrable as well, and
it is hence {\it a fortiori} bounded in $L^2$.
From Lemma~\ref{stm:unfold-en}, restricting the summation in \eqref{eq:elastic} to the set of translations 
\begin{equation*}
	\hat{T}_k \coloneqq \Big\{t\in \inte^3 : \mathrm{dist}\Big(\eps (t+Q), \partial\Omega \Big) > \lambda \eps \Big\} \subset T_k,
\end{equation*}
we deduce
\begin{equation*}
\liminf_{k \to +\infty} \io \chi^0_k(x) W^0_k \big( \eps_k\nv_k(x) P_k^{-1}(x) \big)\dd x 
\geq \liminf_{k \to +\infty} \int_{\hat{\Omega}_k}\int_{Q^0}
W^0_k \big( \nabla_z \hat{v}_k(x,z)P^{-1}_k(x) \big) \dd z \dd x,
\end{equation*}
where
\begin{gather}\label{eq:OmegaQ}
\hat{\Omega}_k \coloneqq \bigcup_{t \in \hat T_k} \eps_k(t+Q).
\end{gather}

We rewrite the right-hand side of the previous inequality as the difference between the integrals
\begin{gather*}\label{diff}
I_k'\coloneqq \int_{\Omega}\int_{Q^0}
W^0_k \big( \nabla_z \hat{v}_k(x,z )P^{-1}_k(x) \big) \dd z \dd x, \\
I_k'' \coloneqq  \int_{\Omega\setminus \hat{\Omega}_k}\int_{Q^0}
W^0_k \big( \nabla_z \hat{v}_k(x,z )P^{-1}_k(x) \big) \dd z \dd x.
\end{gather*}
Being that $\{\nabla_z \hat{v}_k\}$ $2$-equiintegrable, the sequence $\{\nabla_z \hat{v}_k P^{-1}_k\}$ is also $2$-equiintegrable. Thus, by the growth condition \ref{E-growth-eps} and the fact that
$\Omega\setminus \hat{\Omega}_k \subset \{ x \in \Omega : \mathrm{dist}(x,\partial\Omega)\leq (\lambda+\sqrt{3})\eps\} $,
we obtain
$$\lim_{k\to+\infty} I_k'' = 0.$$
Taking into account Lemma~\ref{stm:lem20CC} we conclude
$$
\liminf_{k\to+\infty} \io \chi^0_k(x) W^0_k \big( \eps_k\nv_k(x) P_k^{-1}(x) \big)\dd x \geq 
\liminf_{k\to+\infty} I_k' \geq \mathcal{L}^3(Q^0) \io \rmQ W^0\big(0,P^{-1}(x)\big) \dd x.
$$
\end{proof}


\subsection{Upper bound for the elastic energy}
In this subsection we address the proof of upper $\Gamma$-limit inequality 
for the elastic contribution of the soft component.
Differently from the previous subsection,
in order to establish the desired inequality we perform an analysis that is genuinely two-scale,
in the sense that
we interpret $0$ as the average with respect to the periodic variable
of the two-scale limit of the sequence $\{\eps_k \nv_k\}$.

\begin{proposition}\label{stm:Glimsup-soft}
Let $\{W^0_k\}_{k}$ satisfy assumptions \ref{E-growth-eps}--\ref{E-conv},
and let $P \in W^{1,q}(\Omega;\SL(3))$.
For all $\delta > 0$ there exists a sequence $\{v_k\}\subset W^{1,2}_0(\Omega^0_k;\real^3)$
such that $\eps_k \nabla v_k \rightharpoonup 0$ weakly in $L^2(\Omega;\matr)$
and that
\begin{equation}\label{eq:Glimsup-soft}
\limsup_{k\to +\infty} \io \chi_k^0(x)
W^0_k\left(\eps_k\nabla v_k(x)P^{-1}_k(x)\right) \dd x < \mathcal{L}^3(Q^0) \io \rmQ W^0\big(0,P^{-1}(x)\big) \dd x + \delta,
\end{equation}
whenever $P_k \to P$ uniformly.
\end{proposition}

We begin with a lemma that provides a strong two-scale approximation
of any sufficiently regular function.
The result has already appeared in \cite[Lemma~22]{CC}
where, however, the proof is just sketched.
In Section~\ref{sec:conclusions} we state and prove a more detailed version of this lemma (i.e., Lemma~\ref{stm:recovery-app}) and compare our result with the one in \cite{CC}.

\begin{lemma}\label{stm:recovery}
Let $w\in L^2(\Omega;W^{1,2}_0(Q^0;\real^3)) \cap C^2(\Omega \times Q^0;\real^3)$.
Then, there exists a sequence $\{v_k\} \subset L^2(\Omega;\real^3)$
such that, letting $\hat v_k \coloneqq \mathsf{S}_k v_k$, it holds that
\begin{equation}\label{eq:nablaz-unf}
	\nabla_z \hat v_k \to \nabla_z w
	\quad\text{strongly in } L^2(\Omega\times Q;\matr).
\end{equation}
\end{lemma}

We are now ready to prove the $\Gamma$-limsup inequality
for the soft inclusions functional.

\begin{proof}[Proof of Proposition~\ref{stm:Glimsup-soft}]
According to Lemma~\ref{stm:int-tildeW},
for every $\delta>0$ there exists $w_\delta \in L^2(\Omega;W^{1,2}_0(Q^0;\real^3))$
satisfying
\begin{equation}\label{eq:delta}
\io\int_{Q^0} W^0 \big( \nabla_z w_\delta(x,z) P^{-1}(x) \big) \dd z \dd x
< \mathcal{L}(Q^0) \io \rmQ W^0 \big( 0 , P^{-1}(x) \big) \dd x+\delta
\end{equation}
We would like to apply Lemma~\ref{stm:recovery} which, however, requires $w_\delta \in L^2(\Omega;W^{1,2}_0(Q^0;\real^3))\cap C^2(\Omega \times Q^0;\real^3)$.
We therefore establish the inequality first for a sufficiently regular $w_\delta$,
and we then extend the result by a density argument.

\smallskip
\noindent \textsc{\underline{Case~1: $w_\delta$ regular}}

Let $w_\delta \in L^2(\Omega;W^{1,2}_0(Q^0;\real^3))\cap C^2(\Omega \times Q^0;\real^3)$.
We consider the recovery sequence $\{v_k\}$ coming from Lemma~\ref{stm:recovery}.
Lemmas~\ref{stm:unfolding} and~\ref{lemma:2-scale}(2) yield
$\eps_k \nabla v_k \rightharpoonup 0$ weakly in $L^2(\Omega;\matr)$.
Assumption \ref{E-lip-eps} and H\"older's inequality entail
\begin{multline*}
\sum_{t \in T_k}
\int_{\eps_k(t+Q)} \int_{Q^0} \left\lvert W^0_k \left(
\nabla_z \hat{v}_k (x,z) {P^{-1}_k} (x)
\right) - W^0_k \left(
\nabla_z w_\delta (x,z) {P^{-1}_k} (x)
\right) \right\rvert \dd z \dd x \\
\le c \sum_{t \in T_k}
\left(\int_{\eps_k(t+Q)} \int_{Q^0} \left\lvert
\nabla_z \hat{v}_k (x,z) - 
\nabla_z w_\delta (x,z)
\right\rvert^2 \dd z \dd x\right)^{1/2},
\end{multline*}
where the constant $c$ bounds $\| P^{-1}_k \|_{L^\infty}$.
Thanks to the strong convergence of $\{\nabla_z \hat v_k\}$,
we obtain that the term above is infinitesimal when $k \to +\infty$.
From Lemma~\ref{stm:unfold-en} we then deduce
\begin{multline*}
\limsup_{k\to+\infty} \io \chi_k^0(x)
W^0_k\left(\eps_k\nabla v_k(x)P^{-1}_k(x)\right)\! \dd x
\\	= 	\limsup_{k\to+\infty} \sum_{t \in T_k}
\int_{\eps_k(t+Q)} \int_{Q^0} W^0_k \left(
\nabla_z w_\delta (x,z) {P^{-1}_k} (x)
\right)\! \dd z \dd x
\\	= 	\limsup_{k\to+\infty} \sum_{t \in T_k}
\int_{\eps_k(t+Q)} \int_{Q^0} W^0_k \left(
\nabla_z w_\delta (x,z) {P^{-1}} (x)
\right)\! \dd z \dd x
\\	= 	\limsup_{k\to+\infty} \sum_{t \in T_k}
\int_{\eps_k(t+Q)} \int_{Q^0} W^0 \left(
\nabla_z w_\delta (x,z) {P^{-1}} (x)
\right)\! \dd z \dd x,
\end{multline*}
where the second identity follows from \ref{E-lip-eps} and
the last one from \ref{E-conv}. Note also that, by absolute continuity of the Lebesgue integral,
\begin{multline*}
\limsup_{k\to+\infty} \sum_{t \in T_k}
\int_{\eps_k(t+Q)} \int_{Q^0} W^0 \left(
\nabla_z w_\delta (x,z) {P^{-1}} (x)
\right) \dd z \dd x
\\ = \io \int_{Q^0} W^0 \left(
\nabla_z w_\delta (x,z) {P^{-1}} (x)
\right) \dd z \dd x.
\end{multline*}
Therefore, combining the equalities that we have just found with \eqref{eq:delta}, we achieve the conclusion in the case under consideration.

\smallskip
\noindent \textsc{\underline{Case~2: $w_\delta$ generic}}

\noindent
Let now $w_\delta \in L^2(\Omega;W^{1,2}_0(Q^0;\real^3))$.
By density,
we retrieve a function
$\tilde w_\delta \in C^\infty_c(\Omega; C^\infty_c(Q^0;\real^3))$
such that
\begin{equation*}
\io \int_{Q^0} W^0 \big( \nabla_z \tilde w_\delta(x,z) P^{-1}(x) \big) \dd z
\le \io \int_{Q^0} W^0 \big( \nabla_z w_\delta(x,z) P^{-1}(x) \big) \dd z + \delta.
\end{equation*}
To achieve the conclusion, it only suffices to repeat the argument in Case~1 for $\tilde w_\delta$.
\end{proof}

\subsection{Proof of Proposition~\ref{stm:Glim-soft}}\label{sec:Glimsoft_full}
We are eventually in a position to reap the fruits of the previous subsections.

\begin{proof}[Proof of Proposition~\ref{stm:Glim-soft}]
Let us start with the lower limit inequality.
If the lower limit of $\cJ_k^0(v_k,P_k)$ is not finite,
there is nothing to prove.
Otherwise, statement (1) in Proposition~\ref{stm:Glim-soft} follows
by combining Proposition~\ref{stm:Gliminf-soft} and Lemma~\ref{stm:cont-hard}.

As for the upper bound,
Proposition~\ref{stm:Glimsup-soft} provides for all $\delta > 0$ a sequence $\{v_k\}\subset W^{1,2}_0(\Omega^0_k;\real^3)$
such that
$\eps_k \nabla v_k \rightharpoonup 0$ weakly in $L^2(\Omega;\matr)$
and \eqref{eq:Glimsup-soft} holds.
By the Poincaré inequality on perforated media (see Proposition~\ref{stm:poincare}), 
it follows that $\{v_k\}$ is bounded in $L^2(\Omega;\real^3)$.
We employ again Lemma~\ref{stm:cont-hard} to deduce that
\begin{equation*}
\limsup_{k\to+\infty} \cJ_{k}^0(v_k,P_k) < \cJ^0(P) + \delta.
\end{equation*}
This inequality is actually equivalent to the desired one
(cf.~\cite[Section~1.2]{braides}),
and the proof is therefore concluded.
\end{proof}


\section{Conclusions and a variant of the problem with plastic dissipation}\label{sec:conclusions}
We devote this final section
to the proof of the homogenization result for high-contrast composites
and to the discussion of a variant of the problem
featuring plastic dissipation.

\subsection{Proof of Theorem~\ref{stm:Glim} and convergence of minimum problems}
As we outlined before,
the proof of Theorem~\ref{stm:Glim} is achieved
by combining the splitting procedure in Proposition~\ref{stm:splitting}
with Theorem~\ref{stm:homo-fin-plast2} and Proposition~\ref{stm:Glim-soft},
which account for the asymptotics of the stiff and the soft components, respectively.
Once the homogenization theorem is on hand,
the convergence
of the minimum problems and of their minimizers
will follow thanks to the compactness result in Lemma~\ref{stm:cpt}.

\begin{proof}[Proof of Theorem~\ref{stm:Glim}]
Let $\{\eps_k\}$ be an infinitesimal sequence and
let us fix $y \in  L^2(\Omega;\real^3)$ and $P \in L^q(\Omega;\SL(3))$.
We separate the proof of the lower and of the upper limit inequalities.

\smallskip
\noindent \textsc{\underline{Lower bound}}

\noindent
We consider a sequence $\{(y_k,P_k)\}\subset L^2(\Omega;\real^3) \times L^q(\Omega;\SL(3))$
such that $y_k \to y$ in the sense of extensions and that $P_k \to P$ uniformly.
The only case to discuss is the one in which
the lower limit of $\cJ_k(y_k,P_k)$ is finite, and
we may thus assume that $\{\cJ_k(y_k,P_k)\}$ is bounded.
Keeping in force the notation of Definition~\ref{stm:defconv},
we let $\{\tilde{y}_k\}\subset W^{1,2}(\Omega;\real^3)$ be a sequence such that
$y_k = \tilde y_k$ in $\Omega^1_k$ and
$\tilde y_k \rightharpoonup y$ weakly in $ W^{1,2}(\Omega;\real^3)$.
In the light of \eqref{eq:tildeyk} and Remark~\ref{stm:conv-tilde},
we may without loss of generality assume that
$\tilde y_k \coloneqq \mathsf{T}_k y_k$,
with $\mathsf{T}_k$ as in Lemma~\ref{stm:extension}.

We now apply Proposition~\ref{stm:splitting},
which yields $\{v_k\}\subset W^{1,2}_0(\Omega^0_k;\real^3)$
satisfying \eqref{eq:split-liminf}
and such that $\{v_k\}$ is bounded in $L^2$
and that $\{\eps_k \nv_k\}$ is $2$-equiintegrable.
In particular,
$\eps_k v_k \to 0$ strongly in $L^2$, and hence $(\eps_k v_k,P_k) \stackrel{\tau}{\to} (0,P)$.
Besides, Proposition~\ref{stm:Glim-soft} yields
\begin{equation*}
\cJ^0(P) \leq \liminf_{k\to+\infty} \cJ_k^0(v_k,P_k).
\end{equation*}
At this stage, recalling \eqref{eq:split-liminf},
the proof of the lower bound is concluded
as soon as we show that
\begin{equation}\label{eq:Gliminf-stiff}
\cJ^1(y,P) \leq \liminf_{k\to+\infty} \cJ_k^1(\tilde y_k,P_k) = \liminf_{k\to+\infty} \cJ_k^1(y_k,P_k)
\end{equation}
with $\cJ^1(y,P)$ given by \eqref{eq:J1}.
This is what we prove next.

Let us set
\begin{gather}
\widehat W^1(x,F)\coloneqq \chi_{E^1} (x) W^1(F),
\qquad \widehat H(x,P) \coloneqq \chi_{E^1} (x) H(P), \nonumber \\
\widehat \cJ^1_k(y,P) \coloneqq \io
\left[
\widehat W^1\left(\frac{x}{\eps_k}, \nabla \tilde y P^{-1}\right) 
+
\widehat H \left( \frac{x}{\eps_k},P \right)
+
| \nabla P |^q
\right]\!\dd x. \label{eq:widetilde-Jk}
\end{gather}
It holds
\[
\liminf_{k\to+\infty} \widehat \cJ^1_k(\tilde y_k, P_k)
\leq \liminf_{k\to+\infty} \cJ^1_k(\tilde y_k, P_k).
\]
Since $(\tilde y_k,P_k) \stackrel{\tau}{\to} (y,P)$,
by applying Theorem~~\ref{stm:homo-fin-plast2} to the left-hand side 
of the previous inequality,
\eqref{eq:Gliminf-stiff} is deduced.	 

\smallskip
\noindent \textsc{\underline{Upper bound}}

\noindent
If $P\notin W^{1,q}(\Omega;K)$
there is nothing to prove;
let us then assume that $P \in W^{1,q}(\Omega;K)$.

As we have already observed,
$\{\widehat{\cJ}_k^1\}$ satisfies the requirements of Theorem~\ref{stm:homo-fin-plast2}.
In view of Corollary~\ref{stm:ref-recovery},
for any $(y,P)\in W^{1,2}(\Omega;\real^3) \times W^{1,q}(\Omega;K)$
there exists a sequence $\{(u_k,P_k)\}\subset W^{1,2}(\Omega;\real^3) \times W^{1,q}(\Omega;K)$ such that
$\{\nabla u_k\}$ is $2$-equiintegrable,
$(u_k,P_k) \stackrel{\tau}{\to} (y,P)$, and
$$
\limsup_{k\to+\infty} \widehat{\cJ}_k^1(u_k, P_k)
\leq \cJ^1(y,P).
$$	
Note that
\begin{align*}
0 & \leq \cJ^1_{k}(u_k, P_{k}) - \widehat{\cJ}^1_{k}(u_k, P_{k}) \\
& = \io \big(\chi^1_{k}(x) - \chi_{\eps_k E^1}(x) \big)
\big(
W^1 ( \nabla u_k P^{-1}_{k} ) + H( P_{k} )
\big) \dd x \\
& \leq c \io \big( \chi^1_{k}(x) - \chi_{\eps_k E^1}(x) \big)
\big( |\nabla u_k |^2 + 1 \big) \dd x
\end{align*}
for all $k \in \nat$.
Thanks to the $2$-equiintegrability of $\{\nabla u_k\}$,
we deduce
\begin{equation}\label{eq:limsup-J1}
\limsup_{k\to+\infty} \cJ_k^1(u_k, P_k)
=\limsup_{k\to+\infty} \widehat{\cJ}_k^1(u_k, P_k)
\leq \cJ^1(y,P).
\end{equation}

We now focus on the soft part. Proposition~\ref{stm:Glim-soft} grants the existence 
of a bounded sequence $\{v_k\}\subset L^2(\Omega;\real^3)$ such that $\{v_k\}\subset W^{1,2}_0(\Omega^0_{k};\real^3)$
and
\begin{equation}\label{eq:limsup-J0}
	\limsup_{k\to+\infty} \cJ_{k}^0(v_k,P_k) \leq \cJ^0(P),
\end{equation}
where $\{P_k\}$ is as in \eqref{eq:limsup-J1}.
Notice that if $y_k \coloneqq u_k + v_k$,
then $\{\cJ_k(y_k,P_k)\}$ is bounded and
$\{y_k\}$ converges to $y$ in the sense of extensions
(recall Remark~\ref{stm:conv-tilde}).
Letting $\tilde y_k \coloneqq \mathsf{T}_k y_k$,
thanks to \eqref{eq:split-limsup}
we  conclude the proof of the upper limit inequality:
	\begin{align*}
	\limsup_{k\to+\infty} \cJ_k(y_k,P_k)
			& \leq \limsup_{k\to+\infty} \cJ_k^0(y_k-\tilde y_k,P_k)
			+ \limsup_{k\to+\infty} \cJ_k^1(\tilde y_k,P_k) \\
			& = \limsup_{k\to+\infty} \cJ_k^0(v_k,P_k)
				+ \limsup_{k\to+\infty} \cJ_k^1(u_k,P_k) \\
			& \leq \cJ(y,P).
	\end{align*}
In the previous lines,
the equality is a consequence of the facts that
$\{\nabla u_k\}$ and $\{\nabla \tilde y_k\}$ are bounded and that
$u_k = y_k$ on $\Omega^1_k$,
whereas the last bound is accounted for by \eqref{eq:limsup-J1} and \eqref{eq:limsup-J0}.
\end{proof}

Finally, we are only left to establish the convergence
of the minimum problems associated with the energy functionals $\cJ_\eps$.
What we need is an adaptation of the $\Gamma$-convergence statement
that we have just proved
so as to make it comply with Dirichlet boundary conditions.
To this aim, as it is customary (see, e.g., \cite[Proposition~11.7]{BrDFr}),
we could employ the fundamental estimate derived in \cite{DGP1} on the functionals $\{\widehat{\cJ}^1_k\}$ in \eqref{eq:widetilde-Jk};
indeed, boundary data concern only the stiff part, cf.~Remark~\ref{stm:traces}.
In the light of Corollary~\ref{stm:ref-recovery} we can adopt an alternative strategy.

\begin{proof}[Proof of Corollary~\ref{stm:conv-min}]
Since $\{(y_k,P_k)\}$ is a sequence of almost-minimizers, there exists $C$ such that $\cJ_k(y_k,P_k) \le C$. The $2$-growth condition from below, together with Proposition~\ref{stm:poincare}, provides a bound on $\|y_k\|_{L^2}$. By Lemma~\ref{stm:cpt}, there exists $(y,P) \in W^{1,2}_0(\Omega;\real^3) \times W^{1,q}(\Omega;K)$ such that, up to subsequences, $y_k \to y$ in the sense of extensions and $P_k \to P$ uniformly.
Theorem~\ref{stm:Glim} ensures that
$$\cJ(y,P) \le \liminf_{k\to+\infty} \cJ_k(y_k,P_k).$$

We now prove the existence of a recovery sequence meeting the boundary conditions.
As suggested by Remark~\ref{stm:traces}, we focus on the stiff part.
Let us consider again the functional $\widehat\cJ^1_k$ in \eqref{eq:widetilde-Jk}.
Since the sequence $\{\widehat{\cJ}_k^1\}$ falls within the scopes of Theorem~\ref{stm:homo-fin-plast2},
for any $(\widehat y, \widehat P) \in W^{1,2}_0(\Omega;\real^3) \times W^{1,q}(\Omega;K)$
Corollary~\ref{stm:ref-recovery} provides
a sequence $\{(u_k,\widehat P_k)\}\subset W^{1,2}_0(\Omega;\real^3) \times W^{1,q}(\Omega;K)$ such that
$\{\nabla u_k\}$ is $2$-equiintegrable, $(u_k,\widehat P_k) \stackrel{\tau}{\to} (\widehat y,\widehat P)$ and
$$
	\limsup_{k\to+\infty} \widehat{\cJ}_k^1(u_k,\widehat P_k)
	\leq \cJ^1(y,P).
$$
By reasoning as in the proof of the upper bound in Theorem~\ref{stm:Glim}
we retrieve  a sequence
$\{\widehat y_k,\widehat P_k\}\in W^{1,2}_0(\Omega;\real^3) \times W^{1,q}(\Omega;K)$
such that
$\widehat y_k \to \widehat y$ in the sense of extensions,
$\widehat P_k \to \widehat P$ uniformly and
\begin{align*}
	\limsup_{k\to+\infty} \cJ_k(\widehat y_k,\widehat P_k) \leq \cJ(\widehat y,\widehat P),
\end{align*}
whence
$$
	\limsup_{k\to+\infty} (\inf \cJ_k) \leq \inf \cJ.
$$
Recalling that $\{(y_k,P_k)\}$ is a sequence of almost minimizers, we conclude
$$
	\inf \cJ
	\leq \cJ(y,P)
	\leq \liminf_{k\to+\infty} \cJ_k(y_k,P_k)
	= \liminf_{k \to +\infty} \inf \cJ_k
	\leq  \inf \cJ,
$$
as desired.
\end{proof}

\subsection{A non degenerate upper bound for the soft component}
We proved in Section~\ref{sec:Glimsoft} that
the limiting behavior of the soft inclusions is described
by a degenerate functional.
However, under our assumptions,
a non-degenerate upper bound may still be established,
as we prove in the remainder.
The argument follows \cite{CC},
where {\sc Cherdantsev \& Cherednichenko} derived the effective energy
of high-contrast nonlinear elastic materials.
Differently from us, the $\Gamma$-limit
that they retrieve
keeps track of both the macro- and the microscopic variable, and
this roots in the choice of a stronger notion of convergence.
The drawback of such an approach is the lack of compactness
for sequences with equibounded energy.
It was shown in \cite[Example~2.12]{DKP} that,
when weaker topologies are considered,
the quasiconvex envelope does not provide
the correct limiting energy density for the lower $\Gamma$-limit.

We start by proving a more detailed version of Lemma~\ref{stm:recovery}.
\begin{lemma}[cf.~Lemma~22 in \cite{CC}]\label{stm:recovery-app}
Let $w\in L^2(\Omega;W^{1,2}_0(Q^0;\real^3)) \cap C^2(\Omega \times Q^0;\real^3)$.
Then, there exists a sequence $\{w_k\} \subset L^2(\Omega;W^{1,2}_\per(\real^3;\real^3))$
such that $\nabla_z w_k \to \nabla_z w$ strongly in $L^2(\Omega\times Q;\matr)$.
Besides, setting for $x\in \Omega$
\begin{equation}\label{eq:def-vktilde}
	v_k(x) \coloneqq w_k\bigg( x , \frac{x}{\eps_k}\bigg),
\end{equation}
$\{v_k\}$ converges strongly two-scale to $w$ in $L^2$
and \eqref{eq:nablaz-unf} holds.
\end{lemma}
\begin{proof}
We extend $w$  by setting it equal to $0$ on $Q\setminus Q^0$,
so as to obtain a function in $L^2(\Omega;W^{1,2}_\per(\real^3;\real^3))$ which,
by a slight abuse of notation, we denote again by $w$.

Keeping in mind the definition of $\hat{\Omega}_k$ (see \eqref{eq:OmegaQ}),
for $(\bar x,\bar z) \in \Omega \times \real^3$
we define $w_k(\bar x,\bar z)$ in terms of the averages
of $w(\,\cdot\,,\bar z)$ on the cubes that form $\hat{\Omega}_k$:
\begin{equation}
\label{eq:def-uk}
w_k(\bar x,\bar z) \coloneqq
\begin{cases}
\displaystyle \fint_{\eps_k (t+Q)} w( x , \bar z ) \dd x
& \text{if } \bar x \in \epsilon_k(t+Q) \text{ for some } t \in \hat T_k,\\
0 & \text{for any other } \bar x\in \Omega.
\end{cases}
\end{equation}
By definition, $w_k(\,\cdot\,,z)$ is piecewise constant for all $z\in \bar Q$.
Moreover, for almost every $x\in\Omega$,
$w_k(x,\,\cdot\,)$ is $Q$-periodic as well as weakly differentiable, and
$\nabla_z w_k \to \nabla_z w$ strongly in $L^2(\Omega\times Q;\matr)$.
Indeed, from \eqref{eq:def-uk} and Jensen's inequality, we have that
\begin{align*}
&\io\int_{Q} \left| \nabla_z w_k(x,z) - \nabla_z w(x,z) \right|^2 \dd z \dd x \\
& \quad = \int_{\hat{\Omega}_k} \int_{Q}
\left| \nabla_z w_k(x,z) - \nabla_z w(x,z) \right|^2 
\dd z \dd x
+  \int_{\Omega \setminus \hat{\Omega}_k} \int_{Q}
\left| \nabla_z w(x,z) \right|^2 
\dd z \dd x\\
& \quad = \sum_{t \in \hat{T}_k} \int_{\eps_k(t+Q)} \int_{Q}
\left| \nabla_z w_k(x,z) - \nabla_z w(x,z) \right|^2 
\dd z \dd x
+o(1) \\
& \quad  \leq \sum_{t \in \hat{T}_k} \int_{\eps_k(t+Q)} \int_{Q}
\fint_{\eps_k(t+Q)} \left|  \nabla_z w( \xi , z ) - \nabla_z w\big(x,z) \right|^2 \dd \xi
\dd z \dd x
+o(1),
\end{align*}
and the last term is infinitesimal for $k\to+\infty$
(recall that $w \in C^2$ and the mean value theorem applies). 

We now turn to the functions $v_k$ given by \eqref{eq:def-vktilde}.
First of all we point out that, thanks to the regularity of $w$,
$v_k$ is measurable because it is $C^2$ in the second argument (see \cite[Section~5]{allaire}), and vanishes on $\Omega^1_k$.
Besides, it belongs to $W^{1,2}_0(\Omega^0_k;\real^3)$.
Second, we show that $\{v_k\}$ converges weakly two-scale to $w$ in $L^2$.
To this aim, let us fix $\phi\in C(\bar{\Omega};C_\per(\real^3;\real^3))$.
We find
\begin{align*}
\int_{\Omega} v_k(x)\cdot \phi\left(x,\frac{x}{\eps_k}\right) \dd x 
&  = \int_{\Omega^0_k} w_k\left(x,\frac{x}{\eps_k}\right)\cdot \phi\left(x,\frac{x}{\eps_k}\right) \dd x \\
& = \sum_{t \in T_k}\int_{\eps_k(t+Q^0)} w_k\left(x,\frac{x}{\eps_k}\right)\cdot \phi\left(x,\frac{x}{\eps_k}\right) \dd x \\
& = \eps^3_k\sum_{t \in T_k}\int_{Q^0} w_k\big(\eps_k (t+z),z\big)\cdot \phi\big(\eps_k (t+z),z\big) \dd{z} \\
& = \sum_{t \in \hat{T}_k} \int_{Q^0}\int_{\eps_k(t+Q)} w(x,z)\cdot \phi\big(\eps_k (t+z),z\big) \dd{x}\dd{z} \\
& = \int_{\hat{\Omega}_k} \int_{Q^0} w(x,z) \cdot \phi_k(x,z)\dd z\dd x,
\end{align*}
where $ \phi_k (x,z)\coloneqq \phi (\eps_k(t+z),z)$
if $x\in \eps_k(t+Q)$ with $t\in\hat T_k$.
By the dominated convergence theorem, we infer
\begin{align*}
\lim_{k \to +\infty} \int_{\Omega} v_k(x)\cdot \phi\left(x,\frac{x}{\eps_k}\right) \dd x
= \int_{\Omega}\int_{Q^0} w(x,z) \cdot \phi(x,z) \dd z \dd x,
\end{align*}
that is, $v_k\wts w$ weakly two-scale in $L^2$
(recall that $w(x,z)=0$ if $z\in Q^1$).

In order to prove that strong two-scale convergence actually holds,
we study the limiting behavior of the $L^2$ norm of $\{v_k\}$.
On the one hand, the weak two-scale convergence yields
\begin{align}\label{eq:sci-norme}
\lVert w \rVert_{L^2(\Omega\times Q)}
\leq
\liminf_{k\to+\infty} \lVert v_k \rVert_{L^2(\Omega)}.
\end{align}
On the other hand,
from the properties of $\{w_k\}$ and a change of variables we have the identities
\begin{align*}
\int_{\Omega} | v_k(x)|^2 \dd x 
& = \int_{\Omega^0_k} \left| w_k\left(x,\frac{x}{\eps_k}\right) \right|^2\dd{x}
=\sum_{t\in T_k}\int_{\eps_k(t+Q^0)}\left| w_k\left(x,\frac{x}{\eps_k}\right) \right|^2\dd{x} \\
& =\sum_{t \in T_k} \eps^3_k \int_{Q^0} \left| w_k\big( \eps_k (t+z),z \big) \right|^2\dd{z}
=\sum_{t \in \hat{T}_k} \eps^3_k \int_{Q^0}\left| \fint_{\eps_k(t+Q)}w(x,z)\dd x \right|^2\dd{z}.
\end{align*}
Thanks to Jensen's inequality we deduce
\begin{align*}
\int_{\Omega} \left|v_k(x)\right|^2\dd x \leq
\sum_{t \in \hat{T}_k} \eps^3_k \int_{Q^0}\fint_{\eps_k(t+Q)}|w(x,z)|^2\dd{x}\dd{z}
=\int_{Q^0}\int_{\hat{\Omega}_k}|w(x,z)|^2\dd{x}\dd{z}.
\end{align*}
This, combined with \eqref{eq:sci-norme}, ensures that
$$
\lim_{k\to+\infty}\lVert v_k \rVert_{L^2(\Omega)}
=
\lVert w \rVert_{L^2(\Omega\times Q)}.
$$
In view of Definition~\ref{stm:twoscaleconv}
the conclusion is achieved.

Finally, the strong convergence \eqref{eq:nablaz-unf} follows by observing that, if $x \in \eps_k(t+Q)$, it holds that
\begin{equation*}
\nabla_z \hat v_k(x,z)
= \nabla_z w_k \big( \eps_k (t+z) , z \big).
\end{equation*}
\end{proof}

We are now in a position to prove a non-degenerate upper $\Gamma$-limit inequality
that is the counterpart of the one in Proposition~\ref{stm:Glimsup-soft}
under the current stronger convergence assumptions.

	\begin{proposition}\label{stm:Glimsup-soft-CC}
	Let $\{W^0_k\}_{k}$ satisfy assumptions \ref{E-growth-eps}--\ref{E-conv}.
	For any $(w,P)\in L^2(\Omega;W^{1,2}_0(Q^0;\real^3)) \times W^{1,q}(\Omega;\SL(3))$.
	there exists a sequence $\{v_k\}\subset W^{1,2}_0(\Omega^0_k;\real^3)$
	such that:
	\begin{enumerate}
	\item $v_k \stackrel{2}{\to} w$
		strongly two-scale in $L^2$;
	\item $\eps_k \nv_k \stackrel{2}{\rightharpoonup} \nabla_z w$
		weakly two-scale in $L^2$;
	\item whenever $P_k \to P$ uniformly, it holds that
		\begin{equation*}
			\limsup_{k\to +\infty}
			\io \chi_k^0(x)
				W^0_k \big( \eps_k\nv_k(x)P^{-1}_k(x) \big)
			\dd x
			 \leq \io \int_{Q^0} \rmQ W^0 \big( \nabla_z w(x,z), P^{-1}(x) \big) \dd z \dd x,
		\end{equation*}
		where $\rmQ W^0$ is given by \eqref{eq:tildeW0}. 
	\end{enumerate}
	\end{proposition}

The conclusion is not a straightforward consequence of Lemma~\ref{stm:recovery-app},
because along the sequence $\{v_k\}$ in \eqref{eq:def-vktilde}
we would not end up with the correct limiting energy density.
Therefore, the actual recovery sequence is obtained by adding a ``correction'' to $v_k$.

\begin{proof}[Proof of Proposition~\ref{stm:Glimsup-soft-CC}]
The proof consists of several steps.
At first, to circumvent measurability issues,
it is convenient to consider a sufficiently regular $w$.
Under such assumption,
we are able to construct a recovery sequence of the form
$v_k = \tilde v_k + \tilde w_k,$
where $\{\tilde v_k\}$ is provided by Lemma~\ref{stm:recovery-app}
and $\{\tilde w_k\}$ allows to pass from the densities $W^0_k$ to $\rmQ W^0_k$.
The definition of $\tilde w_k$ is given in Step 1,
while Step 2 deals with the upper limit inequality in the regular case.
The general statement is eventually retrieved by approximation.

\smallskip
\noindent \textsc{\underline{Step 1: construction of $\tilde w_k$ for a regular $w$}}

\noindent
Let us assume that $w\in L^2(\Omega;W^{1,2}_0(Q^0;\real^3))\cap C^2(\Omega \times Q^0;\real^3)$.
We consider a cover of $Q^0$ made of cubes whose edge length is $\eps_k$.
We set $\hat \Sigma_k \coloneqq \set{ s \in \inte^3 : \eps_k (s+Q) \subset \overline{Q^0}}$ and,
for all $(t,s)\in \hat{T}_k \times \hat \Sigma_k$,
we introduce the averages
\begin{equation}\label{eq:def-ak}
A_k(t,s) \coloneqq 
\fint_{\eps_k (t+Q)} \fint_{\eps_k (s+Q)} \nabla_z w( x , z )
\dd z \dd x
\end{equation}
and the piecewise constant functions
\begin{equation*}
A_k(x,z) \coloneqq \begin{cases}
A_k(t,s) &\text{if } (x,z) \in \eps_k(t+Q)\times\eps_k(s+Q),
(t,s)\in \hat{T}_k \times \hat \Sigma_k,\\[3pt]
0	& \text{otherwise}.
\end{cases}
\end{equation*}
We record here for later use that,
by means of Lebesgue differentiation and dominated convergence theorems, it follows
\begin{equation}\label{eq:lim-ak}
\begin{split}
& \lim_{k\to+\infty} \lVert A_k - \nabla_z w 
\rVert^2_{L^2(\Omega\times Q)} \\
& \quad =
\lim_{k\to+\infty} \sum_{t \in \hat{T}_k} \sum_{s \in \hat{\Sigma}_k}
\int_{\eps_k (t+Q)}\int_{\eps_k (s+Q)}
\left| A_k (t,s) - \nabla_z w( x , z )
\right|^2
\dd z \dd x \\
& \quad= 0.
\end{split}
\end{equation}	
By the definition of $\rmQ W^0_k$, for all $k\in \nat$
there exists  $\psi_k \in W^{1,2}_0(Q;\real^3)$ such that
\begin{equation}\label{eq:by-def-tildeW}
\int_{Q} \chi^0(z) W^0_k \Big( \big( A_k(t,s) + \nabla\psi_k(z)\big)  P^{-1}_k( x ) \Big) \dd z \leq \rmQ W^0_k \big( A_k(t,s) , P^{-1}( x )\big) + \frac{1}{k}.
\end{equation}
Note that, due to the smoothness of $w$,
the averages $A_k$ are bounded uniformly in $k$, $t$ and $s$.
In the light of Lemma~\ref{stm:lip},
the values $\rmQ W^0_k \big( A_k(t,s) , P^{-1}( x )\big)$ are uniformly bounded as well.
Therefore,
by combining \eqref{eq:by-def-tildeW}
with assumption \ref{E-growth-eps}, 
we deduce that $\{\psi_k\}$ is bounded in $W^{1,2}_0(Q;\real^3)$.

A change of variables in \eqref{eq:by-def-tildeW} yields
\begin{multline}\label{eq:QpW0k}
\int_{\eps_k(s+Q)} \chi^0 \left( \frac{z}{\eps_k} - s \right) W^0_k \left( \left( A_k(t,s) + \nabla\psi_k \left( \frac{z}{\eps_k} - s \right) \right)  P^{-1}( x ) \right) \dd z  \\ \leq \eps_k^3 \left(\rmQ W^0_k \big( A_k(t,s) , P^{-1}( x )\big) + \frac{1}{k}\right),
\end{multline}
and that suggests us to introduce the functions
\begin{equation*}
\tilde \psi_k(x,z) \coloneqq
\begin{cases}
\eps_k \psi_k \left(\dfrac{z}{\eps_k}-s\right)
&\text{if } (x,z) \in \eps_k(t+Q)\times\eps_k(s+Q), \
(t,s)\in \hat{T}_k \times \hat \Sigma_k,\\[3pt]
0	& \text{otherwise}.
\end{cases}
\end{equation*}
Note that, for each $k$ and $x\in\Omega$,
$\tilde \psi_k(x,\,\cdot\,)$ admits a weak derivative with respect to $z$;
thus, by summing over $(t,s)\in \hat{T}_k \times \hat \Sigma_k$,
from \eqref{eq:QpW0k} we may write
\begin{multline}\label{eq:QpW0k-bis}
\sum_{(t,s)\in \hat{T}_k \times \hat \Sigma_k}
\int_{\eps_k(t+Q)} \int_{\eps_k(s+Q)}
\chi^0 \left( \frac{z}{\eps_k} - s \right) W^0_k \Big( \big( A_k(x,z) + \nabla_z \tilde\psi_k ( x, z ) \big)  P^{-1}( x ) \Big) \dd z \dd x \\
\leq \sum_{(t,s)\in \hat{T}_k \times \hat \Sigma_k} \int_{\eps_k(t+Q)}\eps_k^3 \left(\rmQ W^0_k \big( A_k(t,s) , P^{-1}_k( x )\big) + \frac{1}{k}\right) \dd x.
\end{multline}

We also observe that,
since $\{\psi_k\}$ is bounded,
$\tilde \psi_k \to 0$ strongly in $L^2(\Omega\times Q;\real^3)$.
Then, given that $\{\nabla_z \tilde \psi\}$ is bounded $L^2(\Omega\times Q;\matr)$,
it must converge weakly in $L^2$ to $0$.
It follows that, if $w_k$ is as in Lemma~\ref{stm:recovery-app} and
if $(x,z) \in \eps_k(t+Q)\times\eps_k(s+Q)$
with $(t,s)\in \hat{T}_k \times \hat \Sigma_k$,
\begin{equation}\label{eq:nablaz-uk-wk-bis}
	\nabla_z (w_k+\tilde \psi_k) \rightharpoonup \nabla_z w \quad \text{weakly in } L^2(\Omega\times Q;\matr).
\end{equation}
We further notice that
\begin{align*}
\tilde w_k(x) & \coloneqq \tilde \psi_k\left(x,\frac{x}{\eps_k}\right) \\
&	= \sum_{(t,s)\in \hat{T}_k \times \hat \Sigma_k}
\eps_k\psi_k\left(\dfrac{x}{\eps^2_k}-s\right)
\chi_{\eps_k(t+Q)}(x) \chi_{\eps_k(s+Q)}\left(\frac{x}{\eps_k}\right)
\end{align*}
is a measurable function.
A quick application of the definition of weak derivative proves
also that $\tilde w_k$ belongs to $W^{1,2}_0(\Omega^0_k;\real^3)$.

\smallskip
\noindent \textsc{\underline{Step 2: $w$ regular}}

\noindent We now turn to the proof of the limsup inequality
along the sequence $\{v_k\}$ defined as
\begin{equation}\label{eq:vk}
v_k \coloneqq \tilde v_k + \tilde w_k,
\end{equation}
where
	\[
		\tilde v_k(x)\coloneqq w_k\bigg( x , \frac{x}{\eps_k}\bigg)
	\]
with $w_k$ as in Lemma~\ref{stm:recovery-app}, and
where $\{\tilde w_k\}$ was introduced in Step 1.
We have
\begin{align*}
\hat v_k(x,z) \coloneqq \mathsf{S}_k v_k(x,z)
= w_k \left( \eps_k \left\lfloor\frac{x}{\eps_k} \right\rfloor + \eps_k z , z\right)
+ \tilde \psi_k \left( \eps_k \left\lfloor\frac{x}{\eps_k} \right\rfloor + \eps_k z , z\right),
\end{align*}
so that if $(x,z) \in \eps_k(t+Q)\times\eps_k(s+Q)$
\begin{align}\label{eq:nablaz-hatvk}
\nabla_z \hat v_k(x,z)
& = \nabla_z w_k \big( \eps_k (t+z) , z \big)
+ \nabla\psi_k \left(\dfrac{z}{\eps_k}-s\right).
\end{align}
Taking into account \eqref{eq:nablaz-uk-wk-bis}, \eqref{eq:nablaz-hatvk}
and Lemma~\ref{stm:unfolding}(1),
it follows that
$$
\eps_k \nv_k \stackrel{2}{\rightharpoonup} \nabla_z w
\quad\text{weakly two-scale in } L^2.
$$
Recalling Lemma~\ref{stm:unfold-en}, we have that
\begin{align*}
& \limsup_{k\to +\infty} \io \chi^0_k(x) W^0_k \big( \eps_k\nv_k(x) P^{-1}_k(x) \big)\dd x \\
& \quad=
\limsup_{k\to +\infty} \sum_{t \in T_k}
\int_{\eps_k(t+Q)} \int_{Q^0} W^0_k \left(
\nabla_z \hat{v}_k (x,z) {P^{-1}_k} (x)
\right) \dd z \dd x \\
& \quad = \limsup_{k\to+\infty} I_k,
\end{align*}
where
\begin{align*}
I_k \coloneqq \sum_{(t,s)\in \hat{T}_k \times \hat \Sigma_k}
\int_{\eps_k(t+Q)}\int_{\eps_k(s+Q)}
W^0_k \big(
\nabla_z \hat{v}_k(x , z) P^{-1}_k (x)
\big) \dd z \dd x.
\end{align*}
Indeed,
$\hat v_k$ vanishes
if $x\in \Omega \setminus \hat{\Omega}_k$ or if $z\in Q^0 \setminus \cup \{\epsilon_k(s+Q): s\in \hat\Sigma_k\}$,
and the sequence $\{W^0_k(0)\}$ is bounded by virtue of \ref{E-growth-eps}.
Therefore, since the measure of $\Omega \setminus \hat{\Omega}_k$ and of $Q^0 \setminus \cup \{\epsilon_k(s+Q): s \in \hat\Sigma_k\}$
vanish for $k\to+\infty$,
the second equality holds.

%
%
%

Being the value of $\nabla_z \hat{v}_k \left( x , z \right)$ expressed
by formula \eqref{eq:nablaz-hatvk}, we introduce
\begin{align*}
I_k' \coloneqq
\sum_{t,s}
\int_{\eps_k(t+Q)}\int_{\eps_k(s+Q)}
W^0_k \left(\left(
A_k(t,s)
+ \nabla\psi_k \left(\dfrac{z}{\eps_k}-s\right) \right)
P^{-1}_k(x)
\right) \dd z \dd x,
\end{align*}
where the summation runs over $\hat{T}_k \times \hat \Sigma_k$.
By exploiting assumption \ref{E-lip-eps} and H\"older's inequality,
we obtain the estimate
\begin{align*}
\left| I_k - I'_k \right|
\leq c \sum_{t,s} \int_{\eps_k(t+Q)}\int_{\eps_k(s+Q)}
\left| \Big( \nabla_z w_k \big( \eps_k (t+z) , z \big)  - A_k(t,s) \Big) P^{-1}_k (x)  \right|^2
\dd z \dd x.
\end{align*}
In view of Lemma~\ref{stm:recovery-app} and \eqref{eq:lim-ak}
we deduce
\begin{equation}\label{eq:Ik-Ik'}
\lim_{k\to+\infty} 	\left| I_k - I_k' \right|  = 0.
\end{equation}

Next, let us set
\begin{equation*}
I_k''\coloneqq \int_{\hat{\Omega}_k}\int_{Q^0} \rmQ W^0_k \big( A_k(x,z) , P^{-1}_k(x)\big) \dd z \dd x.
\end{equation*}
According to \eqref{eq:QpW0k-bis},
the difference between the integrands of $I_k'$ and $I_k''$ is of order $k^{-1}$:
\begin{equation}\label{eq:Ik'-Ik''}
\lim_{k\to+\infty} \left| I_k' - I_k'' \right|  = 0.
\end{equation}
Finally, we compare $I''_k$ and the limiting functional.
We have
\begin{align*}
& \left| I_k'' - \io \int_{Q^0} \rmQ W^0 \big( \nabla_z w(x,z) , P^{-1}(x) \big) \dd z \dd x \right|
\\ & \quad\leq 
\int_{\hat{\Omega}_k}\int_{Q^0} 
\left\vert
\rmQ W^0_k \big( A_k(x,z) , P^{-1}_k(x)\big) 
-
\rmQ W^0_k \big( \nabla_z w(x,z) , P^{-1}_k(x) \big)
\right\vert\dd z \dd x \\
& \quad\quad
+
\int_{\hat{\Omega}_k}\int_{Q^0} 
\left\vert
\rmQ W^0_k \big( \nabla_z w(x,z) , P^{-1}_k(x) \big)
-
\rmQ W^0_k \big( \nabla_z w(x,z) , P^{-1}(x) \big)
\right\vert\dd z \dd x \\
& \quad\quad
+
\int_{\hat{\Omega}_k}\int_{Q^0} 
\left\vert
\rmQ W^0_k \big( \nabla_z w(x,z) , P^{-1}(x) \big)
-
\rmQ W^0 \big( \nabla_z w(x,z) , P^{-1}(x) \big)
\right\vert\dd z \dd x \\
& \quad\quad
+
\int_{\Omega\setminus\hat{\Omega}_k}\int_{Q^0} \rmQ W^0 \big( \nabla_z w(x,z) , P^{-1}(x) \big) \dd z \dd x.
\end{align*}
All the terms on the right-hand side vanish as $k\to+\infty$.
Indeed, by using the Lipschitz continuity of $\rmQ W^0_k$
(see Lemma~\ref{stm:lip}(1)) and
the uniform bound on $\{P_k\}$,
the first summand is controlled by the norm of $A_k - \nabla_z v $,
which, according to \eqref{eq:lim-ak}, is infinitesimal.
For what concerns the second term,
Lemma~\ref{stm:lip}(2) and the uniform convergence of $\{P_k\}$ imply
that the integrand is infinitesimal for $k\to+\infty$.
The third quantity vanishes because
$\{\rmQ W^0_k\}$ pointwise converges to $\rmQ W^0$
(recall that they are just variants of the quasiconvex envelopes).
Lastly, the fourth summand is negligible since $\mathcal{L}^3(\Omega\setminus \hat{\Omega}_k)$ tends to $0$.

On the whole, taking into account \eqref{eq:Ik-Ik'} and \eqref{eq:Ik'-Ik''},
we conclude 
\begin{equation*}
\lim_{k\to+\infty} I_k = \io \int_{Q^0} \rmQ W^0 \big( \nabla_z w(x,z), P^{-1}(x) \big) \dd z \dd x.
\end{equation*}

\smallskip
\noindent \textsc{\underline{Step 3: $w$ generic}}

\noindent The argument follows the one of Case~2
in the proof of Proposition~\ref{stm:Glimsup-soft}.

\end{proof}

\subsection{A variant with plastic dissipation}\label{sec:diss}
With a view to applying Theorem~\ref{stm:Glim} to time-dependent problems,
it is useful to modify the functionals $\cJ_\eps$
by adding a term that encodes the plastic dissipation mechanism of the system.
Precisely,
we take into account the non-symmetric distance
$D \colon \matr \times \matr \to [0,+\infty]$ in \eqref{eq:diss-dist}
and we define the dissipation between $P_0,P_1\colon \Omega \to \SL(3)$ as
\begin{equation*}
\cD(P_0; P_1) \coloneqq \io D(P_0, P_1) \dd x.
\end{equation*}
From a physical viewpoint, if $P_0,P_1\colon \Omega \to \SL(3)$ are admissible plastic strains,
$\cD(P_0,P_1)$ is interpreted 
as the minimum amount of energy
that is dissipated when the system moves from a plastic configuration to another.	
Then, assuming that
$\bar P\in W^{1,q}(\Omega;\SL(3))$ represents a pre-existent plastic strain of the body,
we set
\begin{equation}\label{eq:Jeps-diss}
\cJ_\eps^{\mathrm{diss}} (y,P)
\coloneqq \cE_\eps (y,P) + \cD(\bar P;P)+\lVert \nabla P \rVert^q_{L^q(\Omega;\real^{3 \times 3 \times 3})}.
\end{equation}
In the same spirit of \eqref{eq:Eeps0} and \eqref{eq:Eeps1},
we distinguish between the dissipation of the inclusions and the one of the matrix,
respectively
\begin{align*}
\cD_\eps^0(\bar P;P) \coloneqq \io \chi_\eps^0(x) D(\bar P,P) \dd x,
\qquad
\cD_\eps^1(\bar P;P) \coloneqq \io \chi_\eps^1(x) D(\bar P,P) \dd x.
\end{align*}

For what concerns the compactness of sequences with equibounded energy,
we notice that
the presence of the dissipation $\cD$ does not affect Lemma~\ref{stm:cpt}:
the same conclusions hold
if the bound on $\cJ_k(y_k,P_k)$ is replaced
by a bound on $\cJ^{\mathrm{diss}}_k(y_k,P_k)$.

Also our $\Gamma$-convergence results
easily extend to the family $\{\cJ^{\mathrm{diss}}_\eps\}$.
The dissipation is indeed a continuous perturbation:

\begin{lemma}
	Let $P,\bar P\in C(\overline \Omega;K)$ be given.
	If $\{P_k\}\subset C(\overline \Omega;K)$ converges uniformly to $P$,
	then
	$$
	\lim_{k\to+\infty} \cD^i_k(\bar{P};P_k ) = \mathcal{L}^3(Q^i) \cD(\bar{P};P )
	\qquad
	\text{for } i=0,1.
	$$
\end{lemma}
\begin{proof}
	We first observe that if $P_k \to P$ pointwise, then
	\begin{equation}\label{eq:D-to-0}
	D\big( P_k(x) , P(x) \big) \to 0,
	\qquad
	D\big( P(x) , P_k(x) \big) \to 0.
	\end{equation}
	To see this, let $\gamma$ be such that for all $(t,F,G) \in [0,1] \times \SL(3) \times \SL(3)$,
	$\gamma(t,F,G)$ is the evaluation at $t$ of the unique minimizing geodesic
	connecting $F$ and $G$, cf.~Lemma~\ref{stm:K-in-H2}.
	Then, by \eqref{eq:diss-dist} and the definition of $\gamma$,
	\begin{align*}
	D\big( P_k(x) , P(x) \big)
	& = \int_0^1
	\Delta\Big( \gamma \big( t , P_k(x) , P(x) \big) , \dot \gamma \big( t , P_k(x) , P(x) \big) \Big)\!
	\dd t \\
	& \leq c \int_0^1
	| \dot \gamma \big( t , P_k(x) , P(x) \big) |
	\dd t,
	\end{align*}
	where the inequality follows from  the definition of $\Delta$ in \eqref{eq:Delta} and \eqref{Pbound}.
	Since $\dot \gamma$ is continuous and bounded,
	by dominated convergence we deduce that
	the last term vanishes as $k\to +\infty$.
	In a similar fashion, we show that
	$D(P , P_k) \to 0$ as well.
	
	As second step, we notice that	
	\begin{equation}\label{eq:cont-D}
	D\big( \bar P(x) , P_k(x) \big) \to D\big( \bar P(x) , P(x) \big).
	\end{equation}
	Indeed, the triangular inequality yields
	$$
	D\big( \bar P(x) , P(x) \big)-D\big( P_k(x) , P(x) \big)
	\leq D\big( \bar P(x) , P_k(x) \big)
	\leq D\big( \bar P(x) , P(x) \big)+D\big( P(x) , P_k(x) \big),
	$$
	and the assertion follows as a consequence of \eqref{eq:D-to-0}.
	
	Finally, we observe that \eqref{eq:cont-D} grants that
	\begin{align*}
	\lim_{k\to+\infty} \cD^i_k(\bar{P};P_k )
	=
	\lim_{k\to+\infty}
	\int_{\Omega} \chi^i_k(x) D\big( \bar P(x) , P(x) \big)\! \dd x,
	\end{align*}
	and the conclusion is achieved by arguing as in Lemma~\ref{stm:cont-hard}.
\end{proof}

\section*{Acknowledgements}
We acknowledge support
from the Austrian Science Fund (FWF) projects
\href{https://doi.org/10.55776/F65}{10.55776/F65}, \href{https://doi.org/10.55776/V662}{10.55776/V662}, \href{https://doi.org/10.55776/Y1292}{10.55776/Y1292}, and \href{https://doi.org/10.55776/P35359}{10.55776/P35359}, as well as
from the FWF-GA\v{C}R project
\href{https://doi.org/10.55776/I4052}{10.55776/I4052} (19-29646L),
and  
from the OeAD-WTZ project CZ04/2019 (M\v{S}MT\v{C}R 8J19AT013).


\end{document}